\documentstyle[12pt]{article}
\begin{document}
\begin{center}
{\bf ON THE QUANTUM LORENTZ GROUP}\\
\vskip2truecm
M. Lagraa\\
\vskip1truecm
Laboratoire de Physique Th\'eorique\\
Universit\'e d'Oran Es-S\'enia, 31100, Alg\'erie\\
and\\
Laboratoire de Physique Th\'eorique et Hautes \'energies\footnote{Laboratoire
associ\'e au Centre National de la Recherche Scientifique - URA D0063},
Batiment 211,\\Universit\'e Paris XI, 94405 ORSAY,France\\
\vskip1truecm
\end{center}
\vskip3truecm
{\bf Abstract :}
The quantum analogues of Pauli matrices are introduced and investigated. From
these matrices and an appropriate trace over spinorial indices we construct 
a quantum Minkowsky metric. In this framwork we show explicitly the 
correspondence between the $SL(2,C)$ and Lorentz quantum groups. The $\cal R$ 
matrices of the quantum Lorentz group are constructed in terms of the
R matrices of $SL(2,C)$ group. These $\cal R$ matrices satisfy
adequate properties as Yang-Baxter equations, Hecke relations and quantum
symmetrization of the metric. It is also shown that the Minkowsky metric
leads to an invariant and central norm.\\
\newpage
\section{Introduction}
The Lorentz group plays a fundamental role in physics. First, it constitutes
the homogeneous part of Poincar\'e group which is intrinsically connected to
the geometry of the space-time and leaves invariant all physical systems
discribed by the special theory of relativity. Second, the different
representations of Lorentz group are field discribing particles which
constitute the physical systems. For these reasons, it is especially interesting
to study the noncommutative version of the Lorentz group.\\
The other reason which makes the study of the quantum Lorentz group interesting
is due to the fact that in quantum field theory based on a classical space-time
and a classical lorentz group there exist difficulties tied to small
space-time distances. One may hope to solve these difficulties by new tools
provided by the noncommutative geometry [1-4].\\
The construction of quantum Lorentz group has been considered by many authors,
either directly [5-6] or in the context of the quantum poincar\'e group [7-8].
Despite intensive efforts, these studies do not completly describe the quantum
Lorentz group.\\
In this paper one constructs quantum Lorentz group out of the quantum $SL(2,C)$
group by showing how all properties of the former can be deduced from the later
which are well known.\\
The paper is organized in the following way. In Sect. 2 we recall the well
known results provided by the bicovariant calculus over $SL(2,C)$ and $SU(2)$
quantum groups. We shall assume that the undotted (conjugate) and dotted
generators of quantum $SL(2,C)$ group satisfy the commutation rules of the
quantum $SU(2)$ group.\\
In Sect. 3 the construction of quantum Lorentz group is carried out of the
quantum $SL(2,C)$ group following the analogue of the homomorphism for the
classical group $SO(1,3)\sim SL(2,C)\backslash Z_{2}$. We shall start by
investigating the quantum analogues of the Pauli matrices from which we construct
an adequate Minkowsky metric. An example of a two parameters deformation is given
and the completeness relations are established. From the properties of the
quantum Pauli matrices and the generators of $SL(2,C)$, we construct the
generators of the quantum Lorentz group. We show that they satisfy the
axiomatic structure of the Hopf algebras and the orthogonality relations.
We also construct the $\cal R$ matrices of the Lorentz group out of those of
$SL(2,C)$ group. These $\cal R$ matrices satisfy the Yang-Baxter equations,
the Hecke relations and exhibit the quantum symetrization properties of the
Minkowsky metric.\\
In Sect. 4 we investigate the properties of the Minkowsky space.
In particular, we show that the Minkowsky metric induces an invariant norm
which commutes with the Hopf algebra $\cal A$ generated by the quantum
$SL(2,C)$ group generators, the undotted and dotted basis (spinors) of the 
bicovariant bimodule over $\cal A$ and the quantum coordinates of the 
Minkowsky space.\\
\section{Bicovariant Calculus On $SL(2,C)$ and $SU(2)$ Quantum Groups}
Before we start to construct the quantum Lorentz group out of the quantum
$SL(2,C)$ group, let us recall some results provided by the bicovariant calculus
over the $SL(2,C)$ and $SU(2)$ quantum groups. Let an unital $\star$-algebra
$\cal A$ generated by $M_{\alpha}^{~\beta} (\alpha, \beta =1,2)$ which
preserves a nondegenerate bilinear form $\varepsilon$\\
\begin{eqnarray*}
\varepsilon_{\alpha \beta}M_{\gamma}^{~\alpha}M_{\delta}^{~\beta} =
\varepsilon_{\gamma \delta}I_{\cal A}~~,~~\varepsilon^{\gamma \delta}M_{\gamma}^{~\alpha}
M_{\delta}^{~\beta}=\varepsilon^{\alpha \beta}I_{\cal A}~~,~~
\varepsilon^{\alpha \gamma}\varepsilon_{\gamma \beta}=
\delta^{\alpha}_{\beta}=\varepsilon_{\beta \gamma}\varepsilon^{\gamma \alpha}
\end{eqnarray*}
which are the unimodularity conditions. The nondegenerate bilinear form 
$\varepsilon$ is considered as a quantum spinor metric, $I_{\cal A}$ being 
the unity of $\cal A$. To preserve these conditions under the
antimultiplicative involution $\star: \cal A \rightarrow \cal A$,
the spinor metric must satisfy the condition
$(\varepsilon_{\alpha \beta})^{\star} = \lambda \varepsilon_{\dot{\beta} \dot{\alpha}}$
with $\lambda \lambda^{\star}= 1$ leading to:
\begin{eqnarray*}
\varepsilon_{\dot{\alpha} \dot{\beta}}M_{\dot{\gamma}}^{~\dot{\alpha}}
M_{\dot{\delta}}^{~\dot{\beta}} = \varepsilon_{\dot{\gamma} \dot{\delta}}
I_{\cal A}~~,~~
\varepsilon^{\dot{\gamma} \dot{\delta}}M_{\dot{\gamma}}^{~\dot{\alpha}}
M_{\dot{\delta}}^{~\dot{\beta}}=\varepsilon^{\dot{\alpha}\dot{\beta}}
I_{\cal A}~~,~~
\varepsilon^{\dot{\alpha} \dot{\gamma}}\varepsilon_{\dot{\gamma}\dot{\beta}}=
\delta^{\dot{\alpha}}_{\dot{\beta}}=\varepsilon_{\dot{\beta} \dot{\gamma}}
\varepsilon^{\dot{\gamma} \dot{\alpha}}
\end{eqnarray*}
where $M_{\dot{\alpha}}^{~\dot{\beta}}= (M_{\alpha}^{~\beta})^{\star}$.
For convenience, we take
$\lambda=1$. $\cal A$ carries a structure of a $\star$-Hopf algebra with a
coaction $\Delta : \cal A \rightarrow \cal A \otimes \cal A$, a counit
$\varepsilon:  \cal A \rightarrow C$ and an antipode
$S : \cal A \rightarrow \cal A$ defined on the generators by
$\Delta (M_{\alpha}^{~\beta}) = M_{\alpha}^{~\gamma} \otimes
M_{\gamma}^{~\beta}$, $\varepsilon(M_{\alpha}^{~\beta})=
\delta_{\alpha}^{\beta}$ and $S(M_{\alpha}^{~\beta}) =
\varepsilon_{\alpha \gamma} M_{\delta}^{~\gamma} \varepsilon^{\delta \beta}$.
On the dotted copy, we have $ (\Delta(M_{\alpha}^{~\beta}))^{\star} =
\Delta((M_{\alpha}^{~\beta})^{\star})=\Delta (M_{\dot{\alpha}}^{~\dot{\beta}}) =
M_{\dot{\alpha}}^{~\dot{\gamma}} \otimes M_{\dot{\gamma}}^{~\dot{\beta}} $,
$\varepsilon(M_{\dot{\alpha}}^{~\dot{\beta}} )= \delta_{\dot{\alpha}}^{\dot{\beta}}$
and $S(M_{\dot{\alpha}}^{~\dot{\beta}} )=\varepsilon_{\dot{\alpha}\dot{\gamma}}
M_{\dot{\delta}}^{~\dot{\gamma}}\varepsilon^{\dot{\delta}\dot{\beta}} $. the involution
$\star$ acts on the antipode as $(S(M_{\alpha}^{~\beta}))^{\star} =
\varepsilon^{\dot{\beta}\dot{\delta}}M_{\dot{\delta}}^{~\dot{\gamma}}
\varepsilon_{\dot{\gamma}\dot{\alpha}}= S^{-1}(M_{\dot{\alpha}}^{~\dot{\beta}})$.\\
It is known [9] that the generators of a such system satify the noncommutativity
relations $R^{\pm \alpha \beta}_{~~\sigma \rho}M_{\gamma}^{~\sigma}
M_{\delta}^{~\rho} = M_{\sigma}^{~\alpha}M_{\rho}^{~\beta}
R^{\pm \sigma \rho}_{~~\gamma \delta}$ where the forms of the R matrices are
given by $R^{\pm \alpha \beta}_{~\gamma \delta} = \delta^{\alpha}_{\gamma}
\delta^{\beta}_{\delta} + a^{\pm 1} \varepsilon^{\alpha \beta}
\varepsilon_{\gamma \delta}$ satifying $R^{\pm \alpha \beta}_{~\sigma \rho}
R^{\mp \sigma \rho}_{~\gamma \delta}=\delta^{\alpha}_{\gamma}
\delta^{\beta}_{\delta}$ with $a + a^{-1} + \varepsilon^{\alpha \beta}
\varepsilon_{\alpha \beta}= 0$, and $a \not= 0$. These R matrices satisfy the
Yang-Baxter equation, the Hecke equations $(R^{\pm} + a^{\pm 2})(R^{\pm} -1)$
and $\varepsilon_{\alpha \beta}R^{\pm \alpha \lambda}_{~\sigma \gamma}
R^{\pm \beta \rho}_{~\lambda \delta} = a^{\mp 1} \varepsilon_{\gamma \delta}
\delta^{\rho}_{\sigma}$.\\
Now, we consider a right-invariant basis $\theta_{\alpha}$ of the bicovariant
bimodule $\Gamma$ over $\cal A$ on which the right coaction acts as
$\Delta_{R}(\theta_{\alpha}) = \theta_{\alpha} \otimes I$,
$\Delta_{R}(\theta_{\dot{\alpha}})= \theta_{\dot{\alpha}} \otimes I$ and the
left coaction acts as\\
\begin{eqnarray}
\Delta_{L}(\theta_{\alpha}) = M_{\alpha}^{~\beta} \otimes \theta_{\beta}~~&,&~~
\Delta_{L}(\theta^{\alpha}) = S(M_{~\beta}^{~\alpha}) \otimes \theta^{\beta}
\nonumber\\
\Delta_{L}(\theta_{\dot{\alpha}}) =M_{\dot{\alpha}}^{~\dot{\beta}} \otimes \theta_{\dot{\beta}}~~&,&~~
\Delta_{L}(\theta^{\dot{\alpha}}) = S^{-1}(M_{\dot{\beta}}^{~\dot{\alpha}})
\otimes \theta^{\dot{\beta}}
\end{eqnarray}
where $(\theta_{\alpha})^{\star} = \theta_{\dot{\alpha}}$ and the spinorial
indices are lowered and raised as $\theta_{\alpha} =
\theta^{\beta} \varepsilon_{\beta \alpha}, \theta^{\alpha} =
\theta_{\beta} \varepsilon^{\beta \alpha}$, $\theta^{\dot{\alpha}} =
\varepsilon^{\dot{\alpha} \dot{\beta}}\theta_{\dot{\beta}}$ and
$\theta_{\dot{\alpha}} =\varepsilon_{\dot{\alpha} \dot{\beta}}
\theta^{\dot{\beta}} $. From the bicovariance properties of the bimodule
${\cal A}-\Gamma$ [4] we can show the existence of functionals $f: {\cal A}
\rightarrow C$ satisfying the following properties
\begin{eqnarray}
\theta_{\alpha}a = (a \star f_{\alpha}^{~\beta})\theta_{\beta}~~&,&~~
\theta^{\alpha}a = (a \star \tilde{f}_{\beta}^{~\alpha})\theta^{\beta}\\
\theta_{\dot{\alpha}}a = (a \star f_{\dot{\alpha}}^{~\dot{\beta}})
\theta_{\dot{\beta}}~~&,&~~
\theta^{\dot{\alpha}}a = (a \star \tilde{f}_{\dot{\beta}}^{~\dot{\alpha}})\theta^{\dot{\beta}}\\
a\theta_{\alpha} =\theta_{\beta}(a \star f_{\alpha}^{~\beta}\circ S)~~&,&~~
a\theta^{\alpha} =\theta^{\beta}(a \star \tilde{f}_{\beta}^{~\alpha}\circ S)\\
a\theta_{\dot{\alpha}} = \theta_{\dot{\beta}}(a \star
f_{\dot{\alpha}}^{~\dot{\beta}}\circ S)~~&,&~~
a\theta^{\dot{\alpha}} = \theta^{\dot{\beta}}(a \star
\tilde{f}_{\dot{\beta}}^{~\dot{\alpha}}\circ S)\\
f_{\alpha}^{~\beta}(ab)= f_{\alpha}^{~\gamma}(a)f_{\gamma}^{~\beta}(b)~~&,&~~
\tilde{f}_{\alpha}^{~\beta}(ab)= \tilde{f}_{\gamma}^{~\beta}(a)
\tilde{f}_{\alpha}^{~\gamma}(b),\\
f_{\dot{\alpha}}^{~\dot{\beta}}(ab)= f_{\dot{\alpha}}^{~\dot{\gamma}}(a)
f_{\dot{\gamma}}^{~\dot{\beta}}(b)~~&,&~~\tilde{f}_{\dot{\alpha}}^{~\dot{\beta}}(ab)=
\tilde{f}_{\dot{\gamma}}^{~\dot{\beta}}(a)\tilde{f}_{\dot{\alpha}}^{~\dot{\gamma}}(b)\\
f_{\alpha}^{~\beta}(I)=\delta^{\beta}_{\alpha}=\tilde{f}_{\alpha}^{~\beta}(I)~~&,&~~
f_{\dot{\alpha}}^{~\dot{\beta}}(I)=\delta^{\dot{\beta}}_{\dot{\alpha}}=
\tilde{f}_{\dot{\alpha}}^{~\dot{\beta}}(I)\\M_{\alpha}^{~\gamma}(f_{\gamma}^{~\beta}\star a)=
(a \star f_{\alpha}^{~\gamma})M_{\gamma}^{~\beta}~~&,&~~
S(M_{\gamma}^{~\alpha})(\tilde{f}_{\beta}^{~\gamma} \star a)=
(a \star \tilde{f}_{\gamma}^{~\alpha})S(M_{\beta}^{~\gamma}),\\
M_{\dot{\alpha}}^{~\dot{\gamma}}(f_{\dot{\gamma}}^{~\dot{\beta}} \star a)=
(a \star f_{\dot{\alpha}}^{~\dot{\gamma}})M_{\dot{\gamma}}^{~\dot{\beta}}~~&,&~~
S^{-1}(M_{\dot{\gamma}}^{~\dot{\alpha}})
(\tilde{f}_{\dot{\beta}}^{~\dot{\gamma}} \star a)=
(a \star \tilde{f}_{\dot{\gamma}}^{~\dot{\alpha}})
S^{-1}(M_{\dot{\beta}}^{~\dot{\gamma}})
\end{eqnarray}
where the convolution product is defined by $a \star f = (f \otimes I)
\Delta(a)$ for any $a \in \cal A$. Setting $a = S(a)$ ($S^{-1}(a)$) into the
right undotted (dotted) relation of (9)((10)), then applying $S^{-1}$($S$) 
and comparing with the corresponding left undotted (dotted) relation of (9)
((10)), we get
\begin{eqnarray}
f_{\alpha}^{~\beta} = \tilde{f}_{\alpha}^{~\beta}\circ S ~,
f_{\dot{\alpha}}^{~\dot{\beta}} = \tilde{f}_{\dot{\alpha}}^{~\dot{\beta}}\circ S^{-1}.
\end{eqnarray}
For the generators of $\cal A$, the left relation (9) gives
$M_{\alpha}^{~\gamma}M_{\rho}^{~\delta}f_{\gamma}^{~\beta}
(M_{\delta}^{~\sigma})=f_{\alpha}^{~\gamma}(M_{\rho}^{~\delta})
M_{\delta}^{~\sigma}M_{\gamma}^{~\beta}$ which shows that there exist two
functionals $f_{\pm \gamma}^{~~\alpha}(M_{\delta}^{~\rho})$ proportional to
the $R^{\pm \rho \alpha}_{~\gamma \delta}$ matrices. Applying these functionals
on both sides of the unimodularity condition, we obtain:
$f_{\pm \gamma}^{~~\alpha}(M_{\delta}^{~\rho}) = a^{\mp \frac{1}{2}}
R^{\pm \rho \alpha}_{~\gamma \delta}$ and
$f_{\pm \gamma}^{~~\alpha}(S(M_{\delta}^{~\rho})) =
a^{\pm \frac{1}{2}}R^{\mp \alpha \gamma}_{~\delta \gamma}$ [10]. The same
procedure can be used for the dotted copy of $\cal A$ generators. Then there
exist two basis $\theta_{\pm\alpha}$ corresponding to the functionals
$f_{\pm\alpha}^{~~\beta}$.\\
Applying the $\star$ involution on both sides of (2) and (3), we get 
respectively
\begin{eqnarray}
(\theta_{\pm\alpha}a)^{\star} = (\theta_{\pm\beta})^{\star}
(f_{\pm\alpha}^{~~\beta}(a_{(1)}))^{\star}a^{\star}_{(2)}=
a^{\star}(\theta_{\pm\alpha})^{\star} ,~~(\theta_{\pm}^{\alpha}a)^{\star} =
(\theta_{\pm}^{\beta})^{\star}(\tilde{f}_{\pm\beta}^{~~\alpha}
(a_{(1)}))^{\star}a^{\star}_{(2)}=a^{\star}(\theta_{\pm}^{\alpha})^{\star},
\nonumber\\
(\theta_{\pm\dot{\alpha}}a)^{\star} = (\theta_{\pm\dot{\beta}})^{\star}
(f_{\pm\dot{\alpha}}^{~~\dot{\beta}}(a_{(1)}))^{\star}a^{\star}_{(2)}=
a^{\star}(\theta_{\pm\dot{\alpha}})^{\star}~and~
(\theta_{\pm}^{\dot{\alpha}}a)^{\star} = (\theta_{\pm}^{\dot{\beta}})^{\star}
(\tilde{f}_{\pm\dot{\beta}}^{~~\dot{\alpha}}(a_{(1)}))^{\star}a^{\star}_{(2)}=
a^{\star}(\theta_{\pm}^{\dot{\alpha}})^{\star}
\end{eqnarray}
where $a_{(1)}$ and $a_{(2)}$ denote elements of 
$\Delta(a)=a_{(1)} \otimes a_{(2)}$. In the other hand, for $a=a^{\star}$,
(4) and (5) give respectively
\begin{eqnarray}
a^{\star}\theta_{\pm\alpha} = \theta_{\pm\beta}f_{\pm\alpha}^{~~\beta}
(S(a_{(1)}^{\star}))a^{\star}_{(2)}~~,~~
a^{\star}\theta_{\pm}^{\alpha} = \theta_{\pm}^{\beta}
\tilde{f}_{\pm\beta}^{~~\alpha}(S(a_{(1)}^{\star}))a^{\star}_{(2)}\nonumber\\
a^{\star}\theta_{\pm\dot{\alpha}}=\theta_{\pm\dot{\beta}}
f_{\pm\dot{\alpha}}^{~~\dot{\beta}}(S(a_{(1)}^{\star}))a^{\star}_{(2)}~~,~~
a^{\star}\theta_{\pm}^{\dot{\alpha}} = \theta_{\pm}^{\dot{\beta}}
\tilde{f}_{\pm\dot{\beta}}^{~~\dot{\alpha}}(S(a_{(1)}^{\star}))
a^{\star}_{(2)}.
\end{eqnarray}
But for $a = M_{\sigma}^{~\rho}$, we have
$ (f_{\pm\alpha}^{~~\beta}(M_{\sigma}^{~\delta}))^{\star}=
(a^{\mp\frac{1}{2}}R^{\pm \delta\beta}_{~\alpha\sigma})^{\star}=
a^{\mp\frac{1}{2}}R^{\pm \dot{\beta}\dot{\delta}}_{~\dot{\sigma}\dot{\alpha}}=
f_{\mp\dot{\alpha}}^{~~\dot{\beta}}(S(M_{\dot{\sigma}}^{~\dot{\delta}}))$
for $a$ real. Therefore (12) and (13) are consistent if
$(\theta_{\pm\alpha})^{\star} = \theta_{\mp\dot{\alpha}}$ and
$(\theta_{\pm}^{\alpha})^{\star} = \theta_{\mp}^{\dot{\alpha}}$ yielding
\begin{eqnarray}
(f_{\pm\alpha}^{~~\beta}(a))^{\star}=f_{\mp\dot{\alpha}}^{~~\dot{\beta}}
(S(a^{\star})),~(\tilde{f}_{\pm\alpha}^{~~\beta}(a))^{\star}=
\tilde{f}_{\mp\dot{\alpha}}^{~~\dot{\beta}}(S(a^{\star})),\nonumber\\
(f_{\pm\dot{\alpha}}^{~~\dot{\beta}}(a))^{\star}=f_{\mp\alpha}^{~~\beta}
(S(a^{\star}))~and~(\tilde{f}_{\pm\dot{\alpha}}^{~~\dot{\beta}}(a))^{\star}=
\tilde{f}_{\mp\alpha}^{~~\beta}(S(a^{\star})).
\end{eqnarray}
Finally, by using the spinor metric to raise the indices of the right invariant
basis into (2) and (3), we may also show that
\begin{eqnarray}
f_{\pm\alpha}^{~~\beta} = \varepsilon^{\beta\delta}
\tilde{f}_{\pm\delta}^{~~\gamma}\varepsilon_{\gamma\alpha}~~and~~
\varepsilon^{\dot{\beta}\dot{\delta}}f_{\pm\dot{\delta}}^{~~\dot{\gamma}}
\varepsilon_{\dot{\gamma}\dot{\alpha}}=
\tilde{f}_{\pm\dot{\alpha}}^{~~\dot{\beta}}.
\end{eqnarray}
In this stage, we have no indication on the explicit forms of
$f_{\pm \gamma}^{~\alpha}(M_{\dot{\delta}}^{~\dot{\rho}})$ or
$f_{\pm \dot{\gamma}}^{~\dot{\alpha}}(M_{\delta}^{~{\rho}})$ to control
the noncommutativity between undotted and dotted generators of the quantum
$SL(2,C)$ group. To carry this point we assume either the generators
$M_{\alpha}^{~\beta}$ commute with $M_{\dot{\alpha}}^{~\dot{\beta}}$ or are
controled by the R matrices satisfying the properties of the quantum
$SU(2)$ group. In the following we assume the later possibility.\\
To reflect the specific properties of the quantum $SU(2)$ group, we have to add
the unitarity condition on the generators as $M_{\dot{\alpha}}^{~\dot{\beta}} =
S(M_{\beta}^{~\alpha})$ and $M_{\alpha}^{~\beta}=
S^{-1}(M_{\dot{\beta}}^{~\dot{\alpha}})$. Applying $S$ on both sides of
the unimodularity condition, we may show that the unitarity condition yields
$\varepsilon_{\alpha\beta}=\lambda\varepsilon^{\dot{\beta}\dot{\alpha}}$
with $\lambda \lambda^{\star}=1$. In the following we take $\lambda = -1$.
Then to be consistent with the quantum $SU(2)$ group, the spinor metric must
satisfy
\begin{eqnarray}
(\varepsilon_{\alpha \beta})^{\star}=\varepsilon_{\dot{\beta} \dot{\alpha}}=
-\varepsilon^{\alpha \beta}~~and~~
(\varepsilon^{\alpha \beta})^{\star}=\varepsilon^{\dot{\beta} \dot{\alpha}}=
-\varepsilon_{\alpha \beta}.
\end{eqnarray}
It is easy to see, by using the unitarity condition of the generators into (1),
that $\theta_{\dot{\alpha}} = \theta^{\alpha}$ and $\theta_{\alpha}
= \theta^{\dot{\alpha}}$ implying conditions on the functionals 
\begin{eqnarray}
f_{\pm \dot{\alpha}}^{~~\dot{\beta}} = \tilde{f}_{\pm\beta}^{~~\alpha}~~,~~
\tilde{f}_{\pm \dot{\alpha}}^{~~\dot{\beta}} = f_{\pm\beta}^{~~\alpha}
\end{eqnarray}
As stated above, we assume that the functionals of the quantum $SL(2,C)$ group
satisfy the same properties as the $SU(2)$ ones. Therefore, if for example
we set $a= M_{\dot{\sigma}}^{~\dot{\rho}}$ into the left relation of (2), 
we get
\begin{eqnarray}
M_{\alpha}^{~\gamma}M_{\dot{\sigma}}^{~\dot{\delta}}f_{\pm\gamma}^{~~\beta}
(M_{\dot{\delta}}^{~\dot{\rho}}) =f_{\pm \alpha}^{~~\gamma}
(M_{\dot{\sigma}}^{~\dot{\delta}})M_{\dot{\delta}}^{~\dot{\rho}}
M_{\gamma}^{~\beta}~~or~~M_{\alpha}^{~\gamma}M_{\dot{\sigma}}^{~\dot{\delta}}
R^{\pm\dot{\rho}\beta}_{~~\gamma \dot{\delta}}
=R^{\pm \dot{\delta}\gamma}_{~~\alpha\dot{\sigma}}
M_{\dot{\delta}}^{~\dot{\rho}}M_{\gamma}^{~\beta}
\end{eqnarray}
where $f_{\pm \alpha}^{~~\gamma}(M_{\dot{\sigma}}^{~\dot{\delta}})=
R^{\pm\dot{\delta}\gamma}_{~~\alpha\dot{\sigma}} =
f_{\pm \alpha}^{~~\gamma}(S(M_{\delta}^{~\sigma}))=
a^{\pm \frac{1}{2}}R^{\mp \gamma\sigma}_{~~\delta\alpha}$. we have also
objects of the form 
\begin{eqnarray*}
f_{\pm \alpha}^{~~\gamma}
(S(M_{\dot{\sigma}}^{~\dot{\delta}}))=
R^{\mp \gamma\dot{\delta}}_{~~\dot{\sigma}\alpha} =f_{\pm \alpha}^{~~\gamma}
(S(S(M_{\delta}^{~\sigma})))=\varepsilon_{\delta\lambda}
f_{\pm\alpha}^{~~\gamma}
(S(M_{\nu}^{~\lambda}))\varepsilon^{\nu\sigma}= \\
a^{\pm \frac{1}{2}}\varepsilon_{\delta \lambda}R^{\mp \gamma\lambda}
_{~\nu\alpha} \varepsilon^{\nu \sigma}=
\varepsilon_{\dot{\sigma}\dot{\rho}}
f_{\pm\alpha}^{~~\gamma}(M_{\dot{\gamma}}^{~\dot{\rho}})
\varepsilon^{\dot{\gamma}\dot{\delta}}=\varepsilon_{\dot{\sigma}\dot{\rho}}
R^{\pm\dot{\rho}\gamma}_{~~\alpha\dot{\lambda}}
\varepsilon^{\dot{\lambda}\dot{\delta}}. 
\end{eqnarray*}
From (6) and (7), we get respectively
\begin{eqnarray}
f_{\pm\alpha}^{~~\gamma}(M_{\dot{\sigma}}^{~\dot{\delta}})
f_{\pm\gamma}^{~~\beta}(S(M_{\dot{\delta}}^{~\dot{\rho}}))=
\delta_{\alpha}^{\beta}\delta_{\dot{\sigma}}^{\dot{\rho}}=
R^{\pm\dot{\delta}\gamma}_{~~\alpha\dot{\sigma}}
R^{\mp \beta\dot{\rho}}_{~~\dot{\delta}\gamma}\\
f_{\pm\alpha}^{~~\gamma}(S(M_{\dot{\sigma}}^{~\dot{\delta}}))
f_{\pm\gamma}^{~~\beta}(M_{\dot{\delta}}^{~\dot{\rho}})=\delta_{\alpha}^{\beta}
\delta_{\dot{\sigma}}^{\dot{\rho}}=
R^{\mp\gamma\dot{\delta}}_{~~\dot{\sigma}\alpha}
R^{\pm\dot{\rho}\beta}_{~~\gamma\dot{\delta}}
\end{eqnarray}
\section{The quantum Lorentz group}
To have a correspondence between $SL(2,C)$ and Lorentz quantum groups,
we must construct the quantum analogues of the Pauli matrices. Let us consider an
element $X_{\alpha \dot{\beta}}$ as a tensor product of an undotted and
dotted elements of right invariant basis of the bimodule $\cal A$-$\Gamma$
(bispinor). $X_{\alpha \dot{\beta}}$ can always be expanded on a system of
four independent $2 \times 2$ matrices $\sigma^{I}_{~\alpha \dot{\beta}}
(I=0,..,3)$ as $X_{I}\sigma^{I}_{~\alpha \dot{\beta}}$. 
$X_{\alpha \dot{\beta}}$ tranforms as
\begin{eqnarray}
\Delta_{L}(X_{\alpha \dot{\beta}}) = M_{\alpha}^{~\sigma}
M_{\dot{\beta}}^{~\dot{\rho}} \otimes X_{\sigma \dot{\rho}}~~,~~
\Delta_{R}(X_{\alpha \dot{\beta}}) = X_{\alpha \dot{\beta}} \otimes I.
\end{eqnarray}
Then, we have \\
\vskip0.5truecm
{\bf Proposition (3,1)}:\\
$a)~~ There~ exist~ four~ 2~ \times 2~ matrices~
\overline{\sigma}_{\pm}^{I\dot{\alpha} \beta}~ given~ by$
\begin{eqnarray}
\overline{\sigma}_{\pm}^{I\dot{\alpha} \beta} = \varepsilon^{\dot{\alpha}
\dot{\lambda}}R^{\mp \sigma \dot{\rho}}_{~~\dot{\lambda}\nu}
\varepsilon^{\nu \beta}\sigma^{I}_{~\sigma \dot{\rho}}
\end{eqnarray}
$such~ that~ X_{I}\overline{\sigma}_{\pm}^{I\dot{\alpha} \beta} =
X_{\pm}^{~\dot{\alpha} \beta}~ transforms~ under~ the~ left~ coaction~ as$
\begin{eqnarray*}
\Delta_{L}(X_{\pm}^{~\dot{\alpha} \beta}) =
S^{-1}(M_{\dot{\sigma}}^{~\dot{\alpha}})S(M_{\rho}^{~\beta}) \otimes
X_{\pm}^{~\dot{\sigma} \rho}.
\end{eqnarray*}
$b)~~ \overline{\sigma}_{\pm}^{I\dot{\alpha} \beta}~ are~ hermitean~ iff~
\sigma^{I}_{~\alpha\dot{\beta}}~ are.~ In~ this~ case~ X_{I}~ are~ real.$\\
\vskip0.5truecm
$Proof.$: a) Setting into (21) $M_{\alpha}^{~\sigma}
M_{\dot{\beta}}^{~\dot{\rho}} = R^{\pm\dot{\gamma}\delta}_{~\alpha\dot{\beta}}
M_{\dot{\gamma}}^{~\dot{\mu}}M_{\delta}^{~\nu}
R^{\mp\sigma\dot{\rho}}_{~\dot{\mu}\nu}$, obtained by multiplying from the
right both sides of (18) by $R^{\mp\xi\dot{\tau}}_{~\dot{\rho}\beta}$ and
by using (19), we obtain
\begin{eqnarray*}
\Delta_{L}(X_{\alpha \dot{\beta}}) =
R^{\pm\dot{\gamma}\delta}_{~\alpha\dot{\beta}}M_{\dot{\gamma}}^{~\dot{\mu}}
M_{\delta}^{~\nu}R^{\mp\sigma\dot{\rho}}_{~\dot{\mu}\nu} \otimes
X_{\sigma \dot{\rho}}.
\end{eqnarray*}
Multiplying from the left both sides by
$R^{\mp\alpha\dot{\beta}}_{~\dot{\lambda}\nu}$ and using (20), we deduce
\begin{eqnarray*}
\Delta_{L}(X_{I}R^{\mp \alpha\dot{\beta}}_{~\dot{\lambda}\nu}
\sigma^{I}_{~\alpha \dot{\beta}}) = M_{\dot{\lambda}}^{~\dot{\mu}}
M_{\nu}^{~\tau} \otimes X_{I}R^{\mp \sigma\dot{\rho}}_{~\dot{\mu}\tau}
\sigma^{I}_{~\sigma\dot{\rho}}=\\
\varepsilon_{\dot{\lambda}\dot{\gamma}}S^{-1}(M_{\dot{\delta}}^{~\dot{\gamma}})
\varepsilon^{\dot{\delta}\dot{\mu}} \varepsilon^{\tau\xi}
S(M_{\xi}^{~\kappa})\varepsilon_{\kappa\nu}\otimes X_{I}
R^{\mp \sigma\dot{\rho}}_{~\dot{\mu}\tau}\sigma^{I}_{~\sigma\dot{\rho}}
\end{eqnarray*}
yielding
\begin{eqnarray*}
\Delta_{L}(X_{I}\varepsilon^{\dot{\alpha}\dot{\lambda}}
R^{\mp\sigma\dot{\rho}}_{~\dot{\lambda}\nu}\varepsilon^{\nu \beta}
\sigma^{I}_{~\sigma\dot{\rho}})  = S^{-1}(M_{\dot{\delta}}^{~\dot{\alpha}})
S(M_{\gamma}^{~\beta}) \otimes X_{I}\varepsilon^{\dot{\delta}\dot{\lambda}}
R^{\mp\sigma\dot{\rho}}_{~\dot{\lambda}\nu}\varepsilon^{\nu\gamma}
\sigma^{I}_{~\sigma \dot{\rho}},
\end{eqnarray*}
which can be written under the form $\Delta_{L}(X_{\pm}^{~\dot{\alpha}\beta}) =
\Delta_{L}(X_{I}\overline{\sigma}_{\pm}^{I\dot{\alpha} \beta}) =
S^{-1}(M_{\dot{\delta}}^{~\dot{\alpha}}) S(M_{\gamma}^{~\beta}) \otimes
X_{\pm}^{~\dot{\delta} \gamma}$ with
\begin{eqnarray*}
\overline{\sigma}_{\pm}^{I\dot{\alpha} \beta}=
\varepsilon^{\dot{\alpha}\dot{\lambda}}
R^{\mp \sigma\dot{\rho}}_{~\dot{\lambda}\nu}\varepsilon^{\nu \beta}
\sigma^{I}_{\sigma \dot{\rho}}
\end{eqnarray*}
from which we obtain $\varepsilon_{\dot{\lambda}\dot{\alpha}}
\overline{\sigma}_{\pm}^{I\dot{\alpha} \beta}
\varepsilon_{\beta\nu}=
R^{\mp \sigma\dot{\rho}}_{~\dot{\lambda}\nu}
\sigma^{I}_{\sigma \dot{\rho}}$. Multiplying from the right both sides by
$R^{\pm\dot{\lambda}\nu}_{~\gamma\dot{\tau}}$ and using (19), we get
\begin{eqnarray}
\sigma^{I}_{\alpha \dot{\beta}}=\varepsilon_{\dot{\lambda}\dot{\gamma}}
R^{\pm \dot{\lambda}\nu}_{~\alpha\dot{\beta}}\varepsilon_{\mu \nu}
\overline{\sigma}_{\pm}^{I\dot{\gamma} \mu}.
\end{eqnarray}
b) Under the conditions (16) on the spinor metric we have
$(R^{\pm\alpha\rho}_{~~\sigma\nu})^{\star}=R^{\pm\sigma\nu}_{~~\alpha\rho}$,
for $a$ real, implying
\begin{eqnarray*}
(R^{\mp \alpha \dot{\beta}}_{~~\dot{\lambda}\nu})^{\star}=
(a^{\pm \frac{1}{2}}\varepsilon_{\beta \rho}
R^{\mp \alpha \rho}_{~~\sigma\nu} \varepsilon^{\sigma \lambda})^{\star}=
a^{\pm \frac{1}{2}}\varepsilon^{\beta \rho}
R^{\mp \sigma\nu}_{~~\alpha\rho} \varepsilon_{\sigma \lambda}.
\end{eqnarray*}
In the other hand, by using (11) and (15), we obtain 
\begin{eqnarray*}
a^{\pm \frac{1}{2}}\varepsilon^{\beta \rho}R^{\mp \sigma\nu}_{~~\alpha\rho}
\varepsilon_{\sigma \lambda} = f_{\mp\alpha}^{~~\nu}(\varepsilon^{\beta\rho}
M_{\rho}^{~\sigma}\varepsilon_{\sigma\lambda}) = f_{\mp\alpha}^{~~\nu}
(S^{-1}(M_{\lambda}^{~\beta}))=\tilde{f}_{\mp\alpha}^{~~\nu}
(M_{\lambda}^{~\beta})=\\
\varepsilon_{\alpha\sigma}f_{\mp\delta}^{~~\sigma}
(M_{\lambda}^{~\beta})\varepsilon^{\delta\nu} =
a^{\pm \frac{1}{2}}\varepsilon_{\alpha \sigma}
R^{\mp \beta\sigma}_{~~\delta\lambda} \varepsilon^{\delta \nu}=
R^{\mp \beta\ \dot{\alpha}}_{~~\dot{\nu} \lambda}
\end{eqnarray*}
yielding $(\overline{\sigma}_{\pm}^{I\dot{\alpha} \beta})^{\star}=
\varepsilon^{\dot{\beta}\dot{\nu}}
R^{\mp \rho \dot{\sigma}}_{~~\dot{\nu} \lambda}\varepsilon^{\lambda \alpha}
(\sigma^{I}_{~\sigma\dot{\rho}})^{\star}$. Therefore, if
$(\sigma^{I}_{~\sigma\dot{\rho}})^{\star} = \sigma^{I}_{~\rho\dot{\sigma}}$
then
\begin{eqnarray*}
(\overline{\sigma}_{\pm}^{I\dot{\alpha} \beta})^{\star}
=\varepsilon^{\dot{\beta}\dot{\nu}}
R^{\mp \rho\dot{\sigma}}_{~~\dot{\nu} \lambda}\varepsilon^{\lambda \alpha}
\sigma^{I}_{~\rho\dot{\sigma}} =
\overline{\sigma}_{\pm}^{I\dot{\beta} \alpha}.
\end{eqnarray*}
The same procedure can be applyied to (23) to show the converse. It is now
easy to see by applying the $\star$ involution on both sides of (21) that
$(X_{\alpha\dot{\beta}})^{\star} = X_{\beta\dot{\alpha}}$ implying
$X_{I}\sigma^{I}_{\beta\dot{\alpha}}=(X_{I}\sigma^{I}_{\alpha\dot{\beta}})^{\star}=
(X_{I})^{\star}\sigma^{I}_{\beta\dot{\alpha}}$ which shows that $X_{I}$ are
real. Q.E.D. \\
\vskip0.5truecm
To define a quantum metric of the space $\cal M$ spanned by $X_{I}$, we have to define
an adequate trace [10-11] over the spinorial indices which makes invariant
this metric under quantum $SL(2,C)$ group.\\
\vskip0.25truecm
{\bf Proposition (3,2):}\\
${\cal M}~ is~ endowed~ with~ a~ metric~ G^{IJ}~ given~ by$
\begin{eqnarray}
G_{\pm}^{~IJ} =\frac{1}{Q}Tr(\sigma^{I} \overline{\sigma}_{\pm}^{J}) =
\frac{1}{Q}\varepsilon^{\alpha \nu}\sigma^{I}_{~\alpha\dot{\beta}}
\overline{\sigma}_{\pm}^{J\dot{\beta}\gamma}\varepsilon_{\gamma\nu}=
\frac{1}{Q}Tr(\overline{\sigma}_{\pm}^{I} \sigma^{J})=
\frac{1}{Q} \varepsilon_{\dot{\nu}\dot{\gamma}}
\overline{\sigma}_{\pm}^{I\dot{\gamma}\alpha}
\sigma^{J}_{~\alpha\dot{\beta}} \varepsilon^{\dot{\nu}\dot{\beta}}
\end{eqnarray}
$such~that$\\
$a)~~ G_{\pm}^{~IJ}X_{I}X_{J}~ are~ invariant~ under~quantum~ SL(2,C)~ group.$\\
$b)~~ G_{\pm}^{~IJ}~ are~ hermitean~if~ the~ matrices~ 
\sigma^{I}_{~\alpha\dot{\beta}}~ are.$\\
\vskip0.25truecm
$Proof$: a) To show the invariance of $G$ under the quantum $SL(2,C)$ group, 
we consider the norm of $X$ as $G_{\pm}^{~IJ}X_{I}X_{J}=
\frac{1}{Q} \varepsilon^{\alpha\nu}X_{\alpha \dot{\beta}}
X_{\pm}^{~\dot{\beta} \gamma}\varepsilon_{\gamma\nu}$
which transforms under $SL(2,C)$ as
\begin{eqnarray*}
\Delta_{L}(G_{\pm}^{~IJ}X_{I}X_{J})=\frac{1}{Q} \varepsilon^{\alpha\nu}
M_{\alpha}^{~\sigma}M_{\dot{\beta}}^{~\dot{\rho}}
S^{-1}(M_{\dot{\delta}}^{~\dot{\beta}})S(M_{\lambda}^{~\gamma})
\varepsilon_{\gamma\nu} \otimes X_{\sigma \dot{\rho}}
X_{\pm}^{~\dot{\delta}\lambda}=\\
\frac{1}{Q}\varepsilon^{\alpha\nu}M_{\alpha}^{~\sigma}
S(M_{\lambda}^{~\gamma})
\varepsilon_{\gamma\nu} \otimes X_{\sigma \dot{\delta}}
X_{\pm}^{~\dot{\delta}\lambda} = \frac{1}{Q}\varepsilon^{\alpha\nu}
M_{\alpha}^{~\sigma}\varepsilon_{\lambda\rho}M_{\mu}^{~\rho}
\varepsilon^{\mu\gamma}\varepsilon_{\gamma\nu} \otimes
X_{\sigma \dot{\delta}}
X_{\pm}^{~\dot{\delta}\lambda}=\\
I \otimes \frac{1}{Q} \varepsilon^{\sigma\rho}X_{\sigma \dot{\delta}}
X_{\pm}^{~\dot{\delta} \lambda} \varepsilon_{\lambda\rho}=
I \otimes \frac{1}{Q} \varepsilon^{\sigma \rho}
\sigma^{I}_{\sigma\dot{\delta}}
\overline{\sigma}_{\pm}^{J\dot{\delta}\lambda}\varepsilon_{\lambda\rho}
X_{I}X_{J} =I \otimes  G_{\pm}^{~IJ}X_{I}X_{J}.
\end{eqnarray*}
The same computation may be applyed to show that
$\varepsilon_{\dot{\nu}\dot{\gamma}}X_{\pm}^{~\dot{\gamma}\alpha}
X_{\alpha\dot{\beta}}\varepsilon^{\dot{\nu}\dot{\beta}}$ is invariant under
quantum $SL(2,C)$ group.\\
Now, using (22) and the form of the R matrices, we obtain
\begin{eqnarray}
G_{\pm}^{~IJ} = \frac{1}{Q}\varepsilon^{\alpha\xi}
\sigma^{I}_{~\alpha\dot{\beta}}\varepsilon^{\dot{\beta}\dot{\lambda}}
R^{\mp\sigma\dot{\rho}}_{~~\dot{\lambda}\nu}\varepsilon^{\nu\gamma}
\sigma^{J}_{~\sigma\dot{\rho}}\varepsilon_{\gamma\xi}=
\frac{1}{Q}a^{\pm\frac{1}{2}}\varepsilon^{\alpha\xi}
\sigma^{I}_{~\alpha\dot{\beta}}\varepsilon^{\dot{\beta}\dot{\lambda}}\varepsilon_{\rho\delta}R^{\mp\sigma\delta}_{~~\mu\nu}\varepsilon^{\mu\lambda}
\varepsilon^{\nu\gamma}\sigma^{J}_{~\sigma\dot{\rho}}\varepsilon_{\gamma \xi}=\nonumber\\
\frac{1}{Q}a^{\pm\frac{1}{2}}\varepsilon^{\alpha \xi}\varepsilon^{\dot{\beta}\dot{\lambda}}
\varepsilon_{\rho\delta}(\delta^{\sigma}_{\mu}\delta^{\delta}_{\xi} +
a^{\mp 1}\varepsilon^{\sigma\delta}\varepsilon_{\mu\xi})
\varepsilon^{\mu\lambda}\sigma^{I}_{~\alpha\dot{\beta}}\sigma^{J}_{~\sigma\dot{\rho}}
=\nonumber\\
-\frac{1}{Q}a^{\pm\frac{1}{2}}\varepsilon^{\alpha\xi}
\varepsilon^{\dot{\beta}\dot{\lambda}}
\varepsilon_{\rho\delta}(\delta^{\sigma}_{\mu}\delta^{\delta}_{\xi} +
a^{\mp 1}\varepsilon^{\sigma\delta}\varepsilon_{\mu\xi})
\varepsilon_{dot{\lambda}\dot{\mu}}\sigma^{I}_{~\alpha\dot{\beta}}
\sigma^{J}_{~\sigma\dot{\rho}} =\nonumber\\
-\frac{1}{Q}a^{\pm\frac{1}{2}}\varepsilon^{\alpha\xi}
\varepsilon_{\rho\delta}(\delta^{\sigma}_{\mu}\delta^{\delta}_{\xi} +
a^{\mp 1}\varepsilon^{\sigma\delta}\varepsilon_{\mu\xi})
\sigma^{I}_{~\alpha\dot{\mu}}\sigma^{J}_{~\sigma\dot{\rho}} =\nonumber\\
\frac{1}{Q}(a^{\pm\frac{1}{2}}\sigma^{I\xi}_{~~\dot{\mu}}
\sigma^{J~\dot{\xi}}_{~\mu}-a^{\mp\frac{1}{2}}
\sigma^{I\xi\dot{\xi}}\sigma^{J\delta\dot{\delta}})
\end{eqnarray}
where we have used (16) in the third and fifth line and the $\sigma^{I}$ indices
are raised and lowered as for the basis of the bicovariant ${\cal A}-\Gamma$
bimodule ($\sigma^{I\alpha}_{~~\dot{\beta}}=\sigma^{I}_{~\rho\dot{\beta}}
\varepsilon^{\rho\alpha}$, $\sigma^{I~\dot{\beta}}_{~\alpha} =
\varepsilon^{\dot{\beta}\dot{\rho}}\sigma^{I}_{~\alpha\dot{\rho}}$ etc,..).
From a similar computation we can show that
$Tr(\overline{\sigma}_{\pm}^{I} \sigma^{J})$ gives the same form (25) for
$G_{\pm}^{~IJ}$.\\
b) if $\sigma^{I}_{~\alpha\dot{\beta}}$ are hermitean, we have from (16) and
the proposition(3,1)
\begin{eqnarray*}
(G_{\pm}^{~IJ})^{\star} =\frac{1}{Q}
(\varepsilon^{\alpha \nu}\sigma^{I}_{\alpha\dot{\beta}}
\overline{\sigma}_{\pm}^{J\dot{\beta}\gamma}\varepsilon_{\gamma\nu})^{\star}=
\frac{1}{Q}
\varepsilon_{\dot{\nu}\dot{\gamma}}\overline{\sigma}_{\pm}^{J\dot{\gamma}\beta}
\sigma^{I}_{~\beta\dot{\alpha}} \varepsilon^{\dot{\nu}\dot{\alpha}}=
G_{\pm}^{~JI}.
\end{eqnarray*}
which shows that the metric $G_{\pm}^{~IJ}$ is hermitean. Q.E.D.\\
\vskip0.25truecm
Note that if $\sigma^{I}_{\alpha\dot{\beta}}$ are hermitean, then as well as
$X_{I}$ as its norm $G_{\pm}^{IJ}X_{I}X_{J}$ are real.\\
Now, we can give an explicit example where we take
$\sigma^{I}$ the usual four matrices which are the $2 \times 2$
identity matrix $\sigma^{0}_{~\alpha\dot{\beta}}$
and the three Pauli matrices $\sigma^{i}_{~\alpha\dot{\beta}} (i=1,2,3)$ as:\\
\begin{eqnarray*}
\sigma^{0}_{\alpha \dot{\beta}} =\left(\begin{array}{cc}
1 & 0\\
0 & 1
\end{array}
\right)~~,~~
\sigma^{1}_{\alpha \dot{\beta}} =\left(\begin{array}{cc}
0 & 1\\1 & 0
\end{array}
\right)~~,~~
\sigma^{2}_{\alpha \dot{\beta}} =\left(\begin{array}{cc}
0 & -i\\
i &  0
\end{array}
\right)~~,~~
\sigma^{3}_{\alpha \dot{\beta}} =\left(\begin{array}{cc}
1 & 0\\
0 &-1
\end{array}
\right).
\end{eqnarray*}
For the spinorial metric satisfying (16), we may take the most general form as
\begin{eqnarray*}
\varepsilon_{\alpha\beta}=d^{-\frac{1}{2}}\left(\begin{array}{cc}
ir& -q^{-\frac{1}{2}}\\q^{\frac{1}{2}}& ir
\end{array}
\right)~~,~~
\varepsilon^{\alpha\beta}=d^{-\frac{1}{2}}\left(\begin{array}{cc}
ir & q^{-\frac{1}{2}}\\
-q^{\frac{1}{2}}&ir
\end{array}
\right)
\end{eqnarray*}
and
\begin{eqnarray*}
\varepsilon_{\dot{\alpha}\dot{\beta}}=d^{-\frac{1}{2}}\left(\begin{array}{cc}
-ir & q^{\frac{1}{2}}\\
-q^{-\frac{1}{2}}& -ir
\end{array}
\right)~~,~~
\varepsilon^{\dot{\alpha}\dot{\beta}}=d^{-\frac{1}{2}}\left(\begin{array}{cc}-ir & -q^{\frac{1}{2}}\\
q^{-\frac{1}{2}}&-ir
\end{array}
\right)
\end{eqnarray*}
where $q$ and $r\not=\pm1$ are real, $d=-r^{2}+1$ and
$Q=a+a^{-1}= d^{-1}(2r +q + q^{-1}) =-\varepsilon_{\alpha\beta}
\varepsilon^{\alpha\beta}$. With this choice, the computer MAPLE program
gives  metrics $G_{\pm}^{~IJ}$ of the form\\
\vskip0.25truecm
$\left(\begin{array}{clcr}
-A_{(1)}& 0 &2A_{(1)}rd^{-1}(\frac{Q_{(\frac{1}{2})}}{Q})
&A_{(1)}d^{-1}Q_{(-)}\\
0 & A_{(2)} &-iA_{(2)}d^{-1}Q_{(-)}
&2iA_{(2)}rd^{-1}(\frac{Q_{(\frac{1}{2})}}{Q})\\
2A_{(1)}rd^{-1}(\frac{Q_{(\frac{1}{2})}}{Q})&iA_{(2)}d^{-1}Q_{(-)}
&A_{(2)}-4A_{(3)}r^{2}d^{-2}(\frac{(Q_{(\frac{1}{2})})^{2}}{Q})&
-2A_{(2)}rd^{-2}Q_{(-)}Q_{(\frac{1}{2})}\\
A_{(1)}d^{-1}Q{(-)}&-2iA_{(2)}rd^{-1}(\frac{Q_{(\frac{1}{2})}}{Q})
&-2A_{(2)}rd^{-2}Q_{(-)}Q_{(\frac{1}{2})}
&A_{(2)}-A_{(3)}d^{-2}(Q_{(-)})^{2}Q
\end{array}
\right)$ \\
\vskip0.25truecm
where $A_{(1)}=a^{\mp\frac{3}{2}}$, $A_{(2)}=a^{\pm\frac{1}{2}}$,
$A_{(3)}=a^{\mp\frac{1}{2}}$, $Q_{(\frac{1}{2})}=q^{\frac{1}{2}}
+q^{-\frac{1}{2}}$ and $Q_{(-)}=\frac{q-q^{-1}}{Q}$. The inverse $G_{\pm IJ}$
is given by\\
\vskip0.25truecm
$\left(\begin{array}{clcr}
-A_{(4)}\frac{Q^{2}}{4}+\frac{d^{-1}}{4}A_{(2)}(Q_{(\frac{1}{2})})^{2}(Q-2)Q
&0&\frac{rd^{-1}}{2}A_{(3)}Q_{(\frac{1}{2})}Q
&\frac{d^{-1}}{4}A_{(3)}Q_{(-)}Q^{2}\\
0&\frac{Q^{2}}{4}A_{(3)}&i\frac{d^{-1}}{4}A_{(3)}Q_{(-)}Q^{2}
&-i\frac{rd^{-1}}{2}A_{(3)}Q_{(\frac{1}{2})}Q\\
\frac{rd^{-1}}{2}A_{(3)}Q_{(\frac{1}{2})}Q&-i\frac{d^{-1}}{4}A_{(3)}Q_{(-)}Q^{2}&
\frac{Q^{2}}{4}A_{(3)} & 0 \\
\frac{d^{-1}}{4}A_{(3)}Q_{(-)}Q^{2}&
i\frac{rd^{-1}}{2}A_{(3)}Q_{(\frac{1}{2})}Q &0 &\frac{Q^{2}}{4}A_{(3)}
\end{array}
\right)$\\
\vskip0.25truecm
where $A_{(4)} = a^{\pm\frac{3}{2}}$. In the classical limit $r = 0$ and
$q=1$, these metrics  reduce to the classical Minkowsky metric with signature
$(-,+,+,+)$.\\
From the computer MAPLE program, we obtain the completeness relations as
\begin{eqnarray}
\sigma^{I}_{\alpha \dot{\beta}}\overline{\sigma}_{I}^{~\dot{\rho}\sigma}=
Q\delta^{\sigma}_{\alpha}\varepsilon_{\dot{\beta}\dot{\delta}}
\varepsilon^{\dot{\rho}\dot{\delta}}~~,~~
\sigma_{\pm I\alpha \dot{\beta}}\overline{\sigma}_{\pm}^{I\dot{\rho}\sigma}=
Q\delta^{\dot{\rho}}_{\dot{\beta}}\varepsilon_{\delta\alpha}
\varepsilon^{\delta\sigma}
\end{eqnarray}
or
\begin{eqnarray}
\sigma^{I~\dot{\beta}}_{\alpha}\overline{\sigma}_{I\dot{\rho}}^{~~\sigma}=
Q\delta^{\sigma}_{\alpha}\delta_{\dot{\rho}}^{\dot{\beta}}~~,~~
\sigma_{\pm I~\dot{\beta}}^{~~\alpha}\overline{\sigma}^{I~\dot{\rho}}_{\pm\sigma}=Q\delta_{\sigma}^{\alpha}\delta^{\dot{\rho}}_{\dot{\beta}}
\end{eqnarray}
where $\sigma_{\pm I\alpha \dot{\beta}} =
G_{\pm IJ}\sigma^{J}_{\alpha\dot{\beta}}$. Note that from the computer MAPLE 
program we can check that $\overline{\sigma}_{J}^{~\dot{\alpha}\beta}=
G_{+IJ}\overline{\sigma}_{+}^{J\dot{\alpha}\beta}= G_{-IJ}
\overline{\sigma}_{-}^{J\dot{\alpha}\beta}= G_{\pm IJ}
\overline{\sigma}_{\pm}^{J\dot{\alpha}\beta}$ and
$\overline{\sigma}_{+}^{J\dot{\alpha}\beta}G_{+ JI}=
\overline{\sigma}_{-}^{J\dot{\alpha}\beta}G_{- JI}$.\\
\vskip0.25truecm
{$Remark (3,1):$}\\
-- The metric $G_{\pm IJ}$ can be written as $G_{\pm IJ}=
G_{\pm IL}G_{\pm JK}G_{\pm}^{KL}=G_{\pm IL}G_{\pm JK}
\frac{1}{Q}Tr(\sigma^{K}\overline{\sigma}_{\pm}^{L})=
\frac{1}{Q}Tr(\sigma_{\pm J}\overline{\sigma}_{I}) =G_{\pm IL}G_{\pm JK}
\frac{1}{Q}Tr(\overline{\sigma}_{\pm}^{K}\sigma^{L})=
\frac{1}{Q}Tr(\overline{\sigma}_{J}\sigma_{\pm I})$. With a way quite
analogous that was used to derive (25), we may show that
\begin{eqnarray}
G_{\pm IJ}= \frac{1}{Q}(a^{\mp\frac{1}{2}}
\overline{\sigma}_{I~\delta}^{~\dot{\beta}}
\overline{\sigma}_{J\dot{\delta}}^{~~\beta}-
a^{\pm\frac{1}{2}}\overline{\sigma}_{I\dot{\delta}\delta}
\overline{\sigma}_{J\dot{\beta}\beta})
\end{eqnarray}
-- The completeness relations (27) may be used to convert a vector to a
bispinor and vice versa
\begin{eqnarray*}
X_{\alpha\dot{\beta}} = X_{I}\sigma^{I}_{\alpha\dot{\beta}}\Leftrightarrow 
X_{I}= \frac{1}{Q}\varepsilon^{\alpha\nu}X_{\alpha\dot{\beta}}
\overline{\sigma}_{\pm}^{J\dot{\beta}\delta}\varepsilon_{\delta\nu}
G_{\pm JI}~~or ~~
X_{I}= \frac{1}{Q}\varepsilon_{\dot{\nu}\dot{\beta}}
\overline{\sigma}_{I}^{~\dot{\beta}\alpha}X_{\alpha\dot{\delta}}
\varepsilon^{\dot{\nu}\dot{\delta}}.
\end{eqnarray*}
\vskip0.25truecm
We are now ready to show how two copies of undotted and dotted generators of
$\cal A$ may be combined to form generators $\Lambda_{L}^{~K}$ of an unital
algebra $\cal L$ corresponding to the quantum Lorentz group.\\
\vskip0.25truecm
{\bf Theorem (3,1):}\\
$The~ generators~ \Lambda_{L}^{~K}~ (L,K ~= ~0,1,2,3)~of~the~Quantum~Lorentz~
group~are~given~by$
\begin{eqnarray}
\Lambda_{L}^{~K} =\frac{1}{Q}\varepsilon_{\dot{\gamma}\dot{\delta}}
\overline{\sigma}_{L}^{~\dot{\delta}\alpha}M_{\alpha}^{~\sigma}\sigma^{K}_{~\sigma\dot{\rho}}M_{\dot{\beta}}^{~\dot{\rho}}
\varepsilon^{\dot{\gamma}\dot{\beta}}=\frac{1}{Q}
\overline{\sigma}_{L\dot{\gamma}}^{~~\alpha}M_{\alpha}^{~\sigma}
\sigma^{K~\dot{\rho}}_{~\sigma}S^{-1}(M_{\dot{\rho}}^{~\dot{\gamma}})
\end{eqnarray}
$or$
\begin{eqnarray}
\Lambda_{L}^{~K} =\frac{1}{Q}\varepsilon^{\alpha\delta}M_{\alpha}^{~\sigma}
\sigma^{K}_{~\sigma\dot{\rho}}M_{\dot{\beta}}^{~\dot{\rho}}
\overline{\sigma}_{\pm}^{N\dot{\beta}\gamma}
\varepsilon_{\gamma\delta}G_{\pm NL}.
\end{eqnarray}
$They~are~ real~and~ satisfy~ the~ axiomatic~ structure~ of~ Hopf~ algebras~with~
the~relations$
\begin{eqnarray}
G_{\pm NM}\Lambda_{L}^{~N}\Lambda_{K}^{~M} = G_{\pm LK}I_{\cal L}~~and~~
G_{\pm}^{~LK}\Lambda_{L}^{~N}\Lambda_{K}^{~M} = G_{\pm}^{~NM}I_{\cal L}
\end{eqnarray}
where $I_{\cal L}=I_{\cal A}$ is the unity of ${\cal L}\subset {\cal A}$.
\vskip0.25truecm
$Proof$: Multiplying from the left both sides of (21) by
$\overline{\sigma}^{N\dot{\delta}\alpha}$ and making the trace over dotted
indices, we get
\begin{eqnarray*}
\Delta_{L}X_{L}Tr(\overline{\sigma}_{\pm}^{N\dot{\delta}\alpha}
\sigma^{L}_{~\alpha\dot{\beta}})=\Delta_{L}X_{L}QG_{\pm}^{NL}=
\varepsilon_{\dot{\gamma}\dot{\delta}}
\overline{\sigma}_{\pm}^{N\dot{\delta}\alpha}M_{\alpha}^{~\sigma}
\sigma^{K}_{~\sigma\dot{\rho}}M_{\dot{\beta}}^{~\dot{\rho}}
\varepsilon^{\dot{\gamma}\dot{\beta}} \otimes X_{K}
\end{eqnarray*}
yielding
\begin{eqnarray}
\Delta_{L}X_{L}=\frac{1}{Q}\varepsilon_{\dot{\gamma}\dot{\delta}}
\overline{\sigma}_{L}^{~\dot{\delta}\alpha}M_{\alpha}^{~\sigma}
\sigma^{K}_{\sigma\dot{\rho}}M_{\dot{\beta}}^{~\dot{\rho}}
\varepsilon^{\dot{\gamma}\dot{\beta}} \otimes X_{K} = \Lambda_{L}^{~K}
\otimes X_{K}.
\end{eqnarray}
We can also multiply from the right both sides of (21) by
$\overline{\sigma}_{\pm}^{N\dot{\beta}\gamma}$ and take trace over undotted
indices to have (30).\\
To show that (29) is equal to (30), we note that
\begin{eqnarray}
\varepsilon^{\alpha\delta}\overline{\sigma}_{\pm}^{N\dot{\beta}\gamma}
\varepsilon_{\gamma\delta}G_{\pm NL}= \frac{1}{Q}
\varepsilon^{\alpha\delta}\overline{\sigma}_{\pm}^{N\dot{\beta}\gamma}
\varepsilon_{\gamma\delta}\varepsilon_{\dot{\nu}\dot{\mu}}
\overline{\sigma}_{L}^{~\dot{\mu}\rho}\sigma_{\pm N\rho\dot{\tau}}
\varepsilon^{\dot{\nu}\dot{\tau}}=\nonumber\\
\varepsilon^{\alpha\delta}\delta^{\dot{\beta}}_{\dot{\tau}}
\varepsilon_{\sigma\rho}\varepsilon^{\sigma\gamma}
\varepsilon_{\gamma\delta}\varepsilon_{\dot{\nu}\dot{\mu}}
\overline{\sigma}_{L}^{~\dot{\mu}\rho}\varepsilon^{\dot{\nu}\dot{\tau}}=
\varepsilon_{\dot{\nu}\dot{\mu}}\overline{\sigma}_{L}^{~\dot{\mu}\alpha}
\varepsilon^{\dot{\nu}\dot{\beta}}
\end{eqnarray}
where we have used $G_{\pm NL} = Tr(\overline{\sigma}_{L}\sigma_{\pm N})$
and the completeness relation (27). Substituting this equality in (30), we
retrieve (29).\\
The reality of the generators is obtained by noticing that
$(\overline{\sigma}_{\pm}^{N\dot{\beta}\gamma}G_{\pm NL})^{\star}=
G_{\pm LN}\overline{\sigma}_{\pm}^{N\dot{\gamma}\beta}=
\overline{\sigma}_{L}^{~\dot{\gamma}\beta}$ from which, we get
\begin{eqnarray*}
(\Lambda_{L}^{~K})^{\star} =(\frac{1}{Q}\varepsilon^{\alpha\delta}M_{\alpha}^{~\sigma}\sigma^{K}_{~\sigma\dot{\rho}}
M_{\dot{\beta}}^{~\dot{\rho}}\overline{\sigma}_{\pm}^{N\dot{\beta}\gamma}
\varepsilon_{\gamma\delta}G_{\pm NL})^{\star}=\frac{1}{Q}
\varepsilon_{\dot{\delta}\dot{\gamma}}
\overline{\sigma}_{L}^{~\dot{\gamma}\beta}M_{\beta}^{~\rho}
\sigma^{K}_{~\rho\dot{\sigma}}M_{\dot{\alpha}}^{~\dot{\sigma}}
\varepsilon^{\dot{\delta}\dot{\alpha}}=\Lambda_{L}^{~K}
\end{eqnarray*}
due to the fact that $\sigma^{I}_{~\sigma\dot{\rho}}$ are hermitean. The Hopf
structure of the algebra generated by $\Lambda_{L}^{~K}$ is given by:\\
a) Acting the coaction on both sides of (29), we obtain
\begin{eqnarray*}
\Delta(\Lambda_{L}^{~K}) =
\frac{1}{Q}\overline{\sigma}_{L\dot{\gamma}}^{~~\alpha}
M_{\alpha}^{~\delta}S^{-1}(M_{\dot{\nu}}^{~\dot{\gamma}})\otimes
M_{\lambda}^{~\sigma}\sigma^{K~\dot{\rho}}_{~\sigma}
S^{-1}(M_{\dot{\rho}}^{~\dot{\mu}})\delta^{\lambda}_{\delta}
\delta^{\dot{\nu}}_{\dot{\mu}}.
\end{eqnarray*}
From the completeness relation (27), we deduce
\begin{eqnarray}
\Delta(\Lambda_{L}^{~K}) =\frac{1}{Q}
\overline{\sigma}_{L\dot{\gamma}}^{~~\alpha}M_{\alpha}^{~\delta}
\sigma^{I~\dot{\nu}}_{~\delta}S^{-1}(M_{\dot{\nu}}^{~\dot{\gamma}})
\otimes
\frac{1}{Q}\overline{\sigma}_{I\dot{\mu}}^{~~\delta}
M_{\lambda}^{~\sigma}\sigma^{K~\dot{\rho}}_{~\sigma}
S^{-1}(M_{\dot{\rho}}^{~\dot{\mu}})=\Lambda_{L}^{~I}\otimes\Lambda_{I}^{~K}.
\end{eqnarray}
b)The counity acts as
\begin{eqnarray*}
\varepsilon(\Lambda_{L}^{~K}) =\frac{1}{Q}
\overline{\sigma}_{L\dot{\gamma}}^{~~~\alpha}
\varepsilon(M_{\alpha}^{~\sigma})\sigma^{K~\dot{\rho}}_{~\sigma}
\varepsilon(S^{-1}(M_{\dot{\rho}}^{~\dot{\gamma}}))=
\frac{1}{Q}\overline{\sigma}_{L\dot{\gamma}}^{~~~\alpha}
\delta_{\alpha}^{~\sigma}\sigma^{K~\dot{\rho}}_{~\sigma}
\delta_{\dot{\rho}}^{~\dot{\gamma}}= \\
\frac{1}{Q}\overline{\sigma}_{L\dot{\gamma}}^{~~~\alpha}
\sigma^{K~\dot{\gamma}}_{~\alpha}=\frac{1}{Q}G_{\pm LN}
Tr(\overline{\sigma}_{\pm}^{N}\sigma^{K})= G_{\pm LN}G_{\pm}^{NK} =
\delta^{K}_{L}.
\end{eqnarray*}
c) Finally, applying the antipode on both sides of (29), we get
\begin{eqnarray*}
S(\Lambda_{L}^{~K}) =\frac{1}{Q}\varepsilon^{\dot{\gamma}\dot{\beta}}
S(M_{\dot{\beta}}^{~\dot{\rho}})\sigma^{K}_{~\sigma\dot{\rho}}
S(M_{\alpha}^{~\sigma})\overline{\sigma}_{L}^{~\dot{\delta}\alpha}
\varepsilon_{\dot{\gamma}\dot{\delta}}= \frac{1}{Q}
M_{\dot{\nu}}^{~\dot{\gamma}}\varepsilon^{\dot{\nu}\dot{\rho}}
\sigma^{K}_{~\sigma\dot{\rho}}\varepsilon_{\alpha\lambda}M_{\mu}^{~\lambda}
\varepsilon^{\mu\sigma}\overline{\sigma}_{L}^{~\dot{\delta}\alpha}
\varepsilon_{\dot{\gamma}\dot{\delta}}.
\end{eqnarray*}
It now follows from $M_{\dot{\nu}}^{~\dot{\gamma}}M_{\mu}^{~\lambda} =
R^{\mp\tau\dot{\xi}}_{~~\dot{\nu}\mu}M_{\tau}^{~\delta}
M_{\dot{\xi}}^{~\dot{\rho}}R^{\pm\dot{\gamma}\lambda}_{\delta \dot{\rho}}$,
obtained by multiplying from the right both sides of (18)
by $R^{\mp\alpha\dot{\sigma}}_{~\dot{\tau}\nu}$ and by using (19), that
\begin{eqnarray}
S(\Lambda_{L}^{~K})= \frac{1}{Q}R^{\mp\tau\dot{\xi}}_{~~\dot{\nu}\mu}M_{\tau}^{~\delta}M_{\dot{\xi}}^{~\dot{\beta}}
R^{\pm\dot{\gamma}\lambda}_{~~\delta \dot{\beta}}
\varepsilon^{\dot{\nu}\dot{\rho}}\sigma^{K}_{~\sigma\dot{\rho}}
\varepsilon_{\alpha\lambda}\varepsilon^{\mu\sigma}
\overline{\sigma}_{L}^{~\dot{\delta}\alpha}
\varepsilon_{\dot{\gamma}\dot{\delta}}=\nonumber\\
\frac{1}{Q}R^{\mp\tau\dot{\xi}}_{~~\dot{\nu}\mu}M_{\tau}^{~\delta}
M_{\dot{\xi}}^{~\dot{\beta}}\varepsilon^{\dot{\nu}\dot{\rho}}
\sigma^{K}_{~\sigma\dot{\rho}}
\varepsilon^{\mu\sigma}\varepsilon_{\dot{\gamma}\dot{\delta}}
R^{\pm\dot{\gamma}\lambda}_{~~\delta \dot{\beta}}
\varepsilon_{\alpha\lambda}\overline{\sigma}_{L}^{~\dot{\delta}\alpha}=\nonumber\\
\frac{1}{Q}\varepsilon^{\mu\sigma}\sigma^{K}_{~\sigma\dot{\rho}}
\varepsilon^{\dot{\nu}\dot{\rho}}
R^{\mp\tau\dot{\xi}}_{~~\dot{\nu}\mu}M_{\tau}^{~\delta}
M_{\dot{\xi}}^{~\dot{\beta}}\sigma_{\pm L\delta\dot{\beta}} =\nonumber\\
\frac{1}{Q} R^{\mp\tau\dot{\xi}}_{~~\dot{\nu}\mu}
\varepsilon^{\mu\sigma}\sigma^{K~\dot{\nu}}_{~\sigma}M_{\tau}^{~\delta}
M_{\dot{\xi}}^{~\dot{\beta}} \sigma_{\pm L\delta\dot{\beta}}
\end{eqnarray}
where we have used (23) and $ G_{\pm LI}\sigma^{I}_{~\alpha\dot{\beta}}=
\sigma_{\pm L\alpha\dot{\beta}}$ to pass from the second line to the third.
In the other hand, we have
\begin{eqnarray*}
\overline{\sigma}_{N}^{~\dot{\alpha}\beta}G_{\pm}^{NK} = \frac{1}{Q}
\overline{\sigma}_{N}^{~\dot{\alpha}\beta}\varepsilon_{\dot{\nu}\dot{\gamma}}
\overline{\sigma}_{\pm}^{N\dot{\gamma}\delta}\sigma^{K}_{~\delta\dot{\tau}}
\varepsilon^{\dot{\nu}\dot{\tau}}=\frac{1}{Q}
\overline{\sigma}_{N}^{~\dot{\alpha}\beta}
\varepsilon_{\dot{\nu}\dot{\gamma}}\varepsilon^{\dot{\gamma}\dot{\lambda}}
R^{\mp\tau\dot{\rho}}_{~~\dot{\lambda}\sigma}\varepsilon^{\sigma\delta}
\sigma^{N}_{~\tau\dot{\rho}}\sigma^{K}_{~\delta\dot{\tau}}
\varepsilon^{\dot{\nu}\dot{\tau}}
\end{eqnarray*}
where we have used (24) and (22). Using the completeness relation (26),
we obtain
\begin{eqnarray*}
\overline{\sigma}_{N}^{~\dot{\alpha}\beta}G_{\pm}^{NK}=
\delta^{\beta}_{\tau}\varepsilon_{\dot{\rho}\dot{\mu}}
\varepsilon^{\dot{\alpha}\dot{\mu}}R^{\mp\tau\dot{\rho}}_{~~\dot{\nu}\sigma}\varepsilon^{\sigma\delta}\sigma^{K}_{\delta\dot{\tau}}
\varepsilon^{\dot{\nu}\dot{\tau}}
\end{eqnarray*}
or
\begin{eqnarray}
\overline{\sigma}_{N\dot{\alpha}}^{~~~\beta}G_{\pm}^{NK}=
\varepsilon_{\dot{\rho}\dot{\alpha}}
R^{\mp\beta\dot{\rho}}_{~\dot{\nu}\sigma}\varepsilon^{\sigma\delta}
\sigma^{K~\dot{\nu}}_{~\delta}.
\end{eqnarray}
Substituting this equality in (35), we get
\begin{eqnarray*}
S(\Lambda_{L}^{~K})= \frac{1}{Q}\overline{\sigma}_{N\dot{\alpha}}^{~~~\tau}
M_{\tau}^{~\delta}\sigma_{L\delta\dot{\beta}}M_{\dot{\rho}}^{~\dot{\beta}}
\varepsilon^{\dot{\alpha}\dot{\rho}} G_{\pm}^{~NK}=\nonumber\\
\frac{1}{Q}G_{\pm LM}(\overline{\sigma}_{N\dot{\alpha}}^{~~~\tau}
M_{\tau}^{~\delta}\sigma^{M}_{~\delta\dot{\beta}}
M_{\dot{\rho}}^{~\dot{\beta}}\varepsilon^{\dot{\alpha}\dot{\rho}})
G_{\pm}^{NK} =
G_{\pm LM}\Lambda_{N}^{~M}G_{\pm}^{NK}
\end{eqnarray*}
from which we deduce the orthogonality conditions (31). (31) may also be
checked directly by replacing the generators $\Lambda_{L}^{~K}$ by their
expression (29). Then $\Lambda_{L}^{~K}$ generate a Lorentz algebra
$\cal L \subset \cal A$ (Hopf algebra whose generators are subject to Lorentz
group conditions (31)).  Q.E.D.\\
\vskip0.25truecm
The noncommutativity between the generators $\Lambda_{L}^{~K}$ and the 
elements of ${\cal A}$ are given by\\
\vskip0.25truecm
{\bf Theorem (3,2):}\\
$There~ exist~ functionals~ F_{\pm L}^{~~K}: {\cal A} \rightarrow C~ given~
by~$
\begin{eqnarray}
F_{\pm L}^{~~K} = \frac{1}{Q}(\tilde{f}_{\mp \dot{\beta}}^{~~~\dot{\alpha}}
\overline{\sigma}_{L\dot{\alpha}}^{~~~\delta} \star f_{\pm\delta}^{~~~\gamma}
\sigma^{K~\dot{\beta}}_{~\gamma})
\end{eqnarray}
$satisfying:$\\
$a)$
\begin{eqnarray}
\Lambda_{L}^{~I}(F_{\pm I}^{~~K} \star a) =
(a \star F_{\pm L}^{~~I})\Lambda_{I}^{~K}
\end{eqnarray}
$b)$
\begin{eqnarray}
(F_{\pm L}^{~~K}(a))^{\star}=F_{\pm L}^{~~K}(S(a^{\star}))
\end{eqnarray}
$c)$
\begin{eqnarray}
F_{\pm L}^{~~K}(ab) &=& F_{\pm L}^{~~I}(a)F_{\pm I}^{~~K}(b)\\
F_{\pm L}^{~~K}(\varepsilon(a))&=&\delta^{K}_{L}\varepsilon(a)\\
(F_{\pm L}^{~~I} \star F_{\pm I}^{~~K}\circ S)(a)=
(F_{\pm L}^{~~I}\circ S \star F_{\pm I}^{~~K})(a)
&=& \delta^{K}_{L}\varepsilon(a)
\end{eqnarray}
$d)$
\begin{eqnarray}
{\cal R}^{\pm NM}_{~KL}G_{\pm}^{KL} = G_{\pm}^{NM}~~,~~
{\cal R}^{\pm NM}_{~KL}G_{\pm NM} = G_{\pm KL}
\end{eqnarray}
$where~ {\cal R}^{\pm NM}_{~KL} = F_{\pm K}^{~~M}(\Lambda_{L}^{~N})~ satisfy~
the~ Yang-Baxter~ equations~ and~ the~Hecke~ relations$
\begin{eqnarray}
({\cal R}^{\pm} + a^{\pm2})({\cal R}^{\pm} + a^{\mp2})({\cal R}^{\pm} - 1)=0
\end{eqnarray}
\vskip0.25truecm
$Proof$: a) $X_{\alpha\dot{\beta}}$ are a right invariant basis of
${\cal A}$-${\cal M}$ bimodule transforming under the ${\cal L}$ algebra as (21).
Then we can follow the Woronowicz formalism of the bicovariant bimodule [4]
to state
\begin{eqnarray}
X_{\alpha\dot{\beta}}a = (a \star (f_{\mp\dot{\beta}}^{~~\dot{\rho}} \star
f_{\pm\alpha}^{~~\sigma}))X_{\sigma\dot{\rho}} =
f_{\mp\dot{\beta}}^{~~\dot{\rho}}(a_{(1)})f_{\pm\alpha}^{~~\sigma}(a_{(2)})
a_{(3)}X_{\sigma\dot{\rho}},\nonumber\\
aX_{\alpha\dot{\beta}}= X_{\sigma\dot{\rho}}(a \star
(f_{\mp\dot{\beta}}^{~~\dot{\rho}} \star f_{\pm\alpha}^{~~\sigma})\circ S) =
X_{\sigma\dot{\rho}}f_{\pm\alpha}^{~~\sigma}(S(a_{(1)}))
f_{\mp\dot{\beta}}^{~~\dot{\rho}}(S(a_{(2)}))a_{(3)}
\end{eqnarray}
for any $a \in \cal A$. The convolution product of two functional is defined
as $(f_{1} \star f_{2})(a) = (f_{1} \otimes f_{2})\Delta(a)$. This
choice of functionals is justified by the fact that if we apply the $\star$
involution on both sides of the first relation of (45) and we use the
second relation with $a=a^{\star}$, we obtain
\begin{eqnarray}
a^{\star}X_{\beta\dot{\alpha}} = X_{\rho\dot{\sigma}}
(f_{\mp\dot{\beta}}^{~~\dot{\rho}}(a_{(1)}))^{\star}(f_{\pm\alpha}^{~~\sigma}(a_{(2)}))^{\star}(a_{(3)})^{\star}=X_{\rho\dot{\sigma}}
f_{\pm\beta}^{~~\rho}(S(a_{(1)}^{\star}))f_{\mp\dot{\alpha}}^{~~\dot{\sigma}}
(S(a_{(2)}^{\star}))(a_{(3)})^{\star}
\end{eqnarray}
which is consistent with (14). Now, the bicovariant bimodule formalism [4] can 
also be applyied to the $\cal A$-$\cal M$ bimodule to get
\begin{eqnarray*}
M_{\alpha}^{~\rho}M_{\dot{\beta}}^{~\dot{\sigma}}
((f_{\mp\dot{\sigma}}^{~~\dot{\delta}} \star f_{\pm\rho}^{~~\gamma})\star a)=
(a \star (f_{\mp\dot{\beta}}^{~~\dot{\sigma}}\star f_{\pm\alpha}^{~~\rho}))
M_{\rho}^{~\gamma}M_{\dot{\sigma}}^{~\dot{\delta}}
\end{eqnarray*}
from which we deduce
\begin{eqnarray*}
M_{\alpha}^{~\sigma}S^{-1}(M_{\dot{\mu}}^{~\dot{\nu}})
(\tilde{f}_{\mp\dot{\beta}}^{~~\dot{\mu}} \star f_{\pm\sigma}^{~~\gamma}\star a)=
(a \star (\tilde{f}_{\mp\dot{\mu}}^{~~\dot{\mu}}\star f_{\pm\alpha}^{~~\sigma})
M_{\sigma}^{~\gamma}S^{-1}(M_{\dot{\beta}}^{~\dot{\mu}})
\end{eqnarray*}
where we have used (15). Note that this equation is a combination of the first
equation of (9) and the
second equation of (10).\\
Multiplying now both sides from the left by
$\overline{\sigma}_{L\dot{\nu}}^{~~\alpha}$ and from right by
$\sigma^{K~\dot{\beta}}_{~\gamma}$ and using the completeness relation (27), 
we obtain
\begin{eqnarray*}
\frac{1}{Q}\overline{\sigma}_{L\dot{\nu}}^{~~\alpha} M_{\alpha}^{~\sigma}
\sigma^{I~\dot{\mu}}_{~\sigma}S^{-1}(M_{\dot{\mu}}^{~\dot{\nu}})((\tilde{f}_{\mp\dot{\beta}}^{~~\dot{\tau}}
\overline{\sigma}_{I\dot{\tau}}^{~~\delta} \star
f_{\pm\delta}^{~~\gamma}\sigma^{K~\dot{\beta}}_{~\gamma}) \star a)= \\
(a \star (\tilde{f}_{\mp\dot{\mu}}^{~~\dot{\nu}}
\overline{\sigma}_{L\dot{\nu}}^{~~\alpha}\star
f_{\pm\alpha}^{~~\sigma}\sigma^{I~\dot{\mu}}_{~\sigma}))
\frac{1}{Q}\overline{\sigma}_{I\dot{\tau}}^{~~\delta} M_{\delta}^{~\gamma}
\sigma^{K~~\dot{\beta}}_{~\gamma}S^{-1}(M_{\dot{\beta}}^{~\dot{\tau}})
\end{eqnarray*}
yielding
\begin{eqnarray*}
\Lambda_{L}^{~I}(F_{\pm I}^{~~K} \star a) =
(a \star F_{\pm L}^{~~I})\Lambda_{I}^{~K}, with~
F_{\pm L}^{~~K} = \frac{1}{Q}(\tilde{f}_{\mp \dot{\beta}}^{~~\dot{\alpha}}
\overline{\sigma}_{L\dot{\alpha}}^{~~\delta} \star f_{\pm\delta}^{~~\gamma}
\sigma^{K~\dot{\beta}}_{~\gamma}).
\end{eqnarray*}
b) Using (14), (33) and (15), we get
\begin{eqnarray*}
(F_{\pm L}^{~~K}(a))^{\star} = \frac{1}{Q}
(G_{\pm LN}\varepsilon_{\dot{\alpha}\dot{\rho}}
\tilde{f}_{\mp\dot{\beta}}^{~~\dot{\alpha}}(a_{(1)})
\overline{\sigma}_{\pm}^{N\dot{\rho}\delta}f_{\pm\delta}^{~~\gamma}
(a_{(2)})\sigma^{K}_{~\gamma\dot{\lambda}}
\varepsilon^{\dot{\beta}\dot{\lambda}})^{\star}=\\
\frac{1}{Q}\varepsilon^{\lambda\beta}
(f_{\pm\delta}^{~~\gamma}(a_{(2)}))^{\star}\sigma^{K}_{~\lambda\dot{\gamma}}
(\tilde{f}_{\mp\dot{\beta}}^{~~\dot{\alpha}}(a_{(1)}))^{\star}
\overline{\sigma}_{\pm}^{N\dot{\delta}\rho}G_{\pm NL}\varepsilon_{\rho\alpha}=\\
\frac{1}{Q}\varepsilon^{\lambda\beta}f_{\mp\dot{\delta}}^{~~\dot{\gamma}}
(S(a_{(2)}^{\star}))\sigma^{K}_{~\lambda\dot{\gamma}}
\tilde{f}_{\pm\beta}^{~~\alpha}(S(a_{(1)}^{\star}))
\varepsilon_{\dot{\nu}\dot{\mu}}\overline{\sigma}_{L}^{~\dot{\mu}\rho}
\varepsilon_{\alpha\rho}\varepsilon^{\dot{\nu}\dot{\delta}}=\\
\frac{1}{Q}\varepsilon^{\dot{\nu}\dot{\delta}}
f_{\mp\dot{\delta}}^{~~\dot{\gamma}}(S(a_{(2)}^{\star}))
\varepsilon_{\dot{\gamma}\dot{\xi}}
\sigma^{K~\dot{\xi}}_{~\lambda}\varepsilon^{\lambda\beta}
\tilde{f}_{\pm\beta}^{~~\alpha}(S(a_{(1)}^{\star}))\varepsilon_{\alpha\rho}
\overline{\sigma}_{L\dot{\nu}}^{~~~\rho}=\\
\frac{1}{Q}\tilde{f}_{\mp\dot{\xi}}^{~~\dot{\nu}}(S(a_{(2)}^{\star}))
\overline{\sigma}_{L\dot{\nu}}^{~~\rho}f_{\pm\rho}^{~~\lambda}
(S(a_{(1)}^{\star}))\sigma^{K~\dot{\xi}}_{~\lambda}
=F_{\pm L}^{~~K}(S(a^{\star})).
\end{eqnarray*}
c) (40) is deduced directly from (6), (7) and the completeness relation (27) as
\begin{eqnarray*}
F_{\pm L}^{~~K}(ab) =\frac{1}{Q}\tilde{f}_{\mp \dot{\beta}}^{~~\dot{\alpha}}
(a_{(1)}b_{(1)})\overline{\sigma}_{L\dot{\alpha}}^{~~\delta}
f_{\pm\delta}^{~~\gamma}(a_{(2)}b_{(2)})\sigma^{K~\dot{\beta}}_{~\gamma}=\\
\frac{1}{Q}\tilde{f}_{\mp \dot{\nu}}^{~~\dot{\alpha}}(a_{(1)})
\tilde{f}_{\mp\dot{\beta}}^{~~\dot{\nu}}(b_{(1)})
\overline{\sigma}_{L\dot{\alpha}}^{~~\delta}f_{\pm\delta}^{~~\mu}(a_{(2)})
f_{\pm\mu}^{~~\gamma}(b_{(2)})\sigma^{K~\dot{\beta}}_{~\gamma}=\\
\frac{1}{Q}\tilde{f}_{\mp \dot{\rho}}^{~~\dot{\alpha}}(a_{(1)})
\tilde{f}_{\mp\dot{\beta}}^{~~\dot{\nu}}(b_{(1)})
\overline{\sigma}_{L\dot{\alpha}}^{~~\delta}f_{\pm\delta}^{~~\mu}(a_{(2)})
f_{\pm\tau}^{~~\gamma}(b_{(2)})\sigma^{K~\dot{\beta}}_{~\gamma}
\delta^{\tau}_{\mu}\delta^{\dot{\rho}}_{\dot{\nu}}=\\
\frac{1}{Q}(\tilde{f}_{\mp \dot{\rho}}^{~~\dot{\alpha}}
\overline{\sigma}_{L\dot{\alpha}}^{~~\delta} \star f_{\pm\delta}^{~~\mu}
\sigma^{I~\dot{\rho}}_{~\mu})(a)
\frac{1}{Q}(\tilde{f}_{\mp \dot{\beta}}^{~~\dot{\nu}}
\overline{\sigma}_{I\dot{\nu}}^{~~\tau} \star f_{\pm\tau}^{~~\gamma}
\sigma^{K~\dot{\beta}}_{~\gamma})(b)=F_{\pm L}^{~~I}(a) F_{\pm I}^{~~I}(b).
\end{eqnarray*}
We also have
\begin{eqnarray*}
F_{\pm L}^{~~K}(\varepsilon (a)) =\frac{1}{Q}
\tilde{f}_{\mp \dot{\beta}}^{~~\dot{\alpha}}(\varepsilon (a_{(1)}))
\overline{\sigma}_{L\dot{\alpha}}^{~~\delta}f_{\pm\delta}^{~~\gamma}
(\varepsilon (a_{(2)}))\sigma^{K~\dot{\beta}}_{~\gamma} =
\frac{1}{Q}\overline{\sigma}_{L\dot{\alpha}}^{~~~\delta}
\sigma^{K~\dot{\alpha}}_{~\delta}\varepsilon(a)=\delta^{K}_{L}
\varepsilon(a)
\end{eqnarray*}
where we have used (8). (42) can be deduced directly from (40) and (41).\\
d) Applying (37) on (29), we get
\begin{eqnarray*}
{\cal R}^{\pm NM}_{~LK}G_{\pm}^{LK}=\frac{1}{Q}
(\tilde{f}_{\mp \dot{\beta}}^{~~\dot{\tau}}
\overline{\sigma}_{L\dot{\tau}}^{~~\delta} \star f_{\pm\delta}^{~~\gamma}
\sigma^{M~\dot{\beta}}_{~\gamma})(\frac{1}{Q}
\overline{\sigma}_{K\dot{\nu}}^{~~~\sigma}M_{\sigma}^{~\rho}
\sigma^{N~\dot{\lambda}}_{~\rho}S^{-1}(M_{\dot{\lambda}}^{~\dot{\nu}}))
G_{\pm}^{~LK}= \\
\frac{1}{Q^{2}} \tilde{f}_{\mp \dot{\beta}}^{~~\dot{\tau}}
(M_{\sigma}^{~\alpha}S^{-1}(M_{\dot{\mu}}^{~\dot{\nu}}))
\overline{\sigma}_{L\dot{\tau}}^{~~~\delta}
\overline{\sigma}_{K\dot{\nu}}^{~~~\sigma}f_{\pm\delta}^{~~\gamma}
(M_{\alpha}^{~\rho}S^{-1}(M_{\dot{\lambda}}^{~\dot{\mu}}))
\sigma^{M~\dot{\beta}}_{~\gamma} \sigma^{N~\dot{\lambda}}_{~\rho}
G_{\pm}^{~LK}.
\end{eqnarray*}
Using $ G_{\pm}^{~LK}\overline{\sigma}_{L\dot{\tau}}^{~~~\delta}
\overline{\sigma}_{K\dot{\nu}}^{~~~\sigma}=
\overline{\sigma}_{L\dot{\tau}}^{~~~\delta}
\overline{\sigma}_{\pm~\dot{\nu}}^{L~~\sigma}=
R^{\pm\alpha\dot{\beta}}_{~\dot{\nu}\mu}\varepsilon^{\mu\sigma}
\sigma^{L}_{~\alpha\dot{\beta}}
\overline{\sigma}_{L\dot{\tau}}^{~~\delta} =-Q(a^{\pm\frac{1}{2}}
\varepsilon^{\delta\nu} \varepsilon^{\tau\sigma} + a^{\mp\frac{1}{2}}
\varepsilon^{\delta\tau} \varepsilon^{\sigma\nu})$, obtained from (22), the
completeness relation (27), (16) and the form of the R matrices, we get 
\begin{eqnarray*}
-\frac{1}{Q} \tilde{f}_{\mp \dot{\kappa}}^{~~\dot{\tau}}(M_{\sigma}^{~\alpha})
\tilde{f}_{\mp \dot{\beta}}^{~~\dot{\kappa}}(S^{-1}(M_{\dot{\mu}}^{~\dot{\nu}}))
f_{\pm\delta}^{~~\xi}(M_{\alpha}^{~\rho})f_{\pm\xi}^{~~\gamma}
(S^{-1}(M_{\dot{\lambda}}^{~\dot{\mu}}))(a^{\pm\frac{1}{2}}
\varepsilon^{\delta\nu}\varepsilon^{\tau\sigma} + a^{\mp\frac{1}{2}}
\varepsilon^{\delta\tau} \varepsilon^{\sigma\nu})=\\
-\frac{1}{Q}R^{\mp\alpha\kappa}_{~\tau\sigma}R^{\mp\mu\beta}_{~\kappa\nu}
R^{\pm\rho\xi}_{~\delta\alpha}R^{\pm\lambda\gamma}_{~\xi\mu}(a^{\pm\frac{1}{2}}\varepsilon^{\delta\nu} \varepsilon^{\tau\sigma} +
a^{\mp\frac{1}{2}}\varepsilon^{\delta\tau}\varepsilon^{\sigma\nu})
\end{eqnarray*}
where we have used (17). An explicit computation gives
\begin{eqnarray*}
R^{\mp\alpha\kappa}_{~\tau\sigma}R^{\mp\mu\beta}_{~\kappa\nu}
R^{\pm\rho\xi}_{~\delta\alpha}R^{\pm\lambda\gamma}_{~\xi\mu}
\varepsilon^{\delta\nu}\varepsilon^{\tau\sigma}=
\varepsilon^{\rho\beta}\varepsilon^{\lambda\gamma}~and~
R^{\mp\alpha\kappa}_{~\tau\sigma}R^{\mp\mu\beta}_{~\kappa\nu}
R^{\pm\rho\xi}_{~\delta\alpha}R^{\pm\lambda\gamma}_{~\xi\mu}
\varepsilon^{\delta\tau} \varepsilon^{\sigma\nu}=\varepsilon^{\rho\lambda}
\varepsilon^{\gamma\beta}
\end{eqnarray*}
from which we deduce
\begin{eqnarray*}
{\cal R}^{\pm NM}_{~LK}G^{\pm LK}=-\frac{1}{Q}(a^{\pm\frac{1}{2}}
\varepsilon^{\rho\beta}\varepsilon^{\lambda\gamma}+a^{\mp\frac{1}{2}}
\varepsilon^{\rho\lambda}\varepsilon^{\gamma\beta})
\sigma^{M~\dot{\beta}}_{~\gamma}\sigma^{N~\dot{\lambda}}_{~\rho}=\\
\frac{1}{Q}(a^{\pm\frac{1}{2}}\varepsilon^{\rho\beta}\varepsilon_{\dot{\gamma}\dot{\lambda}}-a^{\mp\frac{1}{2}}
\varepsilon^{\rho\lambda}\varepsilon^{\gamma\beta})
\sigma^{M~\dot{\beta}}_{~\gamma}\sigma^{N~\dot{\lambda}}_{~\rho}=\\
\frac{1}{Q}(a^{\pm\frac{1}{2}}\sigma^{N\beta}_{~~~\dot{\gamma}}
\sigma^{M~\dot{\beta}}_{~\gamma}-a^{\mp\frac{1}{2}}
\sigma^{N\lambda\dot{\lambda}}\sigma^{M\beta\dot{\beta}})=G_{\pm}^{NM}
\end{eqnarray*}
where we have used (16) and (25). A similar calculation gives
\begin{eqnarray*}
{\cal R}^{\pm NM}_{~LK}G_{\pm NM}=\frac{1}{Q}(a^{\mp\frac{1}{2}}
\overline{\sigma}_{L~\delta}^{~\dot{\beta}}
\overline{\sigma}_{K\dot{\delta}}^{~~~\beta}-a^{\pm\frac{1}{2}}
\overline{\sigma}_{L\dot{\beta}\beta}\overline{\sigma}_{K\dot{\delta}\delta})=
G_{\pm LK}.
\end{eqnarray*}
The Yang-Baxter equations may be obtained by applying $F_{\pm L}^{~~K}$ on
both sides of (38) for $a=\Lambda_{N}^{~M}$ and then using (40). Using the
Hecke relations of the R matrices, we obtain after a explicit straightforward 
calculation 
\begin{eqnarray*}
{\cal R}^{\pm NM}_{~IJ}{\cal R}^{\pm IJ}_{~LK}=(2 - a^{\pm 2}-
a^{\mp 2}){\cal R}^{\pm NM}_{~LK}+\delta_{L}^{N}\delta_{K}^{M}+\\
\frac{1}{Q^{2}}a^{\mp 2}(1-a^{\pm 2})R^{\pm\rho\xi}_{~\delta\alpha}
R^{\pm\lambda\gamma}_{~\xi\kappa}R^{\mp\alpha\kappa}_{~\tau\sigma}
\sigma^{M~\dot{\beta}}_{~\gamma}\sigma^{N~\dot{\lambda}}_{~\rho}
\overline{\sigma}_{L\dot{\tau}}^{~~\delta}
\overline{\sigma}_{K\dot{\beta}}^{~~\sigma}+\\
\frac{1}{Q^{2}}a^{\pm 2}(1-a^{\mp 2})R^{\mp\mu\beta}_{~\xi\nu}
R^{\pm\lambda\gamma}_{~\alpha\mu}R^{\mp\alpha\xi}_{~\tau\sigma}
\sigma^{M~\dot{\beta}}_{~\gamma}\sigma^{N~\dot{\lambda}}_{~\rho}
\overline{\sigma}_{L\dot{\tau}}^{~~\rho}
\overline{\sigma}_{K\dot{\nu}}^{~~\sigma},\\
\frac{1}{Q^{2}}R^{\pm\rho\xi}_{~\delta\alpha}
R^{\pm\lambda\gamma}_{~\xi\kappa}R^{\mp\alpha\kappa}_{~\tau\sigma}
\sigma^{M~\dot{\beta}}_{~\gamma}\sigma^{N~\dot{\lambda}}_{~\rho}
\overline{\sigma}_{I\dot{\tau}}^{~~\delta}
\overline{\sigma}_{J\dot{\beta}}^{~~\sigma}{\cal R}^{\pm IJ}_{~LK}=\\
(1-a^{\pm 2}){\cal R}^{\pm NM}_{~LK}+\frac{a^{\pm 2}}{Q^{2}}
R^{\mp\mu\beta}_{~\xi\nu}R^{\pm\lambda\gamma}_{~\kappa\mu}
R^{\mp\kappa\xi}_{~\tau\sigma}\sigma^{M~\dot{\beta}}_{~\gamma}
\sigma^{N~\dot{\lambda}}_{~\rho}\overline{\sigma}_{L\dot{\tau}}^{~~\rho}
\overline{\sigma}_{K\dot{\nu}}^{~~\sigma}~~and\\
\frac{1}{Q^{2}}R^{\mp\mu\beta}_{~\xi\nu}R^{\pm\lambda\gamma}_{~\alpha\mu}
R^{\mp\alpha\xi}_{~\tau\sigma}\sigma^{M~\dot{\beta}}_{~\gamma}
\sigma^{N~\dot{\lambda}}_{~\rho}\overline{\sigma}_{I\dot{\tau}}^{~~\rho}
\overline{\sigma}_{J\dot{\nu}}^{~~\sigma}{\cal R}^{\pm IJ}_{~LK}=\\
(1-a^{\mp 2}){\cal R}^{\pm NM}_{~LK}+ \frac{a^{\mp 2}}{Q^{2}}
R^{\pm\delta\alpha}_{~\xi\kappa}R^{\pm\lambda\gamma}_{~\alpha\mu}
R^{\mp\kappa\mu}_{~\tau\sigma}\sigma^{M~\dot{\beta}}_{~\gamma}
\sigma^{N~\dot{\delta}}_{~\rho}\overline{\sigma}_{L\dot{\tau}}^{~~\xi}
\overline{\sigma}_{K\dot{\beta}}^{~~\sigma}
\end{eqnarray*}
leading to (44) . Q.E.D.\\
\vskip0.25truecm
As for the ${\cal A}$-$\Gamma$ bimodule, the existence of two functionals 
$F_{\pm L}^{~~K}$ leads to two right invariant basis $X_{\pm I}$ of the 
${\cal A}-{\cal M}$ bimodule satisfying
\begin{eqnarray}
X_{\pm L}a=(a \star F_{\pm L}^{~~K})X_{\pm K}&,&~ X_{\pm}^{(\pm)L}a=
(a \star \tilde{F}_{\pm K}^{(\pm)~L})X_{\pm}^{(\pm)K}\nonumber\\
aX_{\pm L}=X_{\pm K}(a \star F_{\pm L}^{~~K}\circ S)&,&~aX_{\pm}^{(\pm)L}=
X_{\pm}^{(\pm)K}(a \star \tilde{F}_{\pm K}^{(\pm)~L}\circ S)
\end{eqnarray}
for any $a\in \cal A$. The indices are raised and lowered by using the
Minkowsky metric as $X_{L}G_{\pm}^{LK}=X^{(\pm)K}$. Following the same
formalism applyied to the ${\cal A}$-$\Gamma$ bimodule, we may  show that
\begin{eqnarray*}
F_{\pm L}^{~~K} = G_{\pm}^{KM}\tilde{F}_{\pm M}^{(\pm)~N}G_{\pm NL}~and~
\tilde{F}_{\pm L}^{(\pm)~K} \circ S = F_{\pm L}^{~~K}
\end{eqnarray*}
implying
\begin{eqnarray}
G_{\pm}^{~MN}F_{\pm N}^{~~L} \star F_{\pm M}^{~~K}(a) = G_{\pm}^{~KL}
\varepsilon(a)~and~G_{\pm KL}F_{\pm N}^{~~L} \star F_{\pm M}^{~~K}(a) =
G_{\pm MN}\varepsilon(a)
\end{eqnarray}
for any $a\in \cal A$.
\section{\bf Quantum Minkowsky space}
We consider the elements $X_{I}$ of right invariant basis of the 
$\cal A$-$\cal M$ bimodule as coordinates which span the Minkowsky space 
$\cal M$ over the field R. Then Minkowski space is equipped either with 
coordinates $X_{+ I}$ and metric $G_{+ IJ}$ or with coordinates $X_{- I}$ 
and metric $G_{- IJ}$. (32) shows that 
$\Delta_{L}:{\cal M} \rightarrow {\cal L} \otimes {\cal M}$ is a
corepresentation of ${\cal L}$ in the vector space ${\cal M}$. In fact, from 
(32), (34) and the action of the counity on the generators of $\cal L$, it is 
easy to verify $(id \otimes \Delta_{L})\Delta_{L}=(\Delta \otimes id)
\Delta_{L}$ and $(\varepsilon \otimes id)\Delta_{L} = id$.\\
In the following we  assume that the coordinates $X_{\pm I}$ commute , in 
the quantum sens, with them selves and with the elements (spinors) of the right 
invariant basis of $\cal A$-$\Gamma$ bimodule. To carry this quantum 
symmetrization, we consider the bicovariant bimodule automorphism $\sigma $ [4] 
such that for any $a,b\in \cal A$ and the left invariant element
$\tilde{X}_{\pm I} = S(\Lambda_{I}^{~J})X_{\pm J}$,
$\tilde{\theta}_{\pm\alpha} = S(M_{\alpha}^{~\beta})\theta_{\pm\beta}$
or $\tilde{\theta}_{\pm\dot{\alpha}} = S(M_{\dot{\alpha}}^{~\dot{\beta}})
\theta_{\pm\dot{\beta}}$, we have
\begin{eqnarray*}
\sigma (\tilde{X}_{(a) L}\otimes X_{(b) K})=X_{(b)K}\otimes \tilde{X}_{(a)L}=
X_{(b)K}\otimes S(\Lambda_{L}^{~N})  X_{(a) N} =\\
S(\Lambda_{L}^{~N})\sigma (X_{(a) N} \otimes X_{(b) K})=
F_{(b)K}^{~~~Q}(S(\Lambda_{P}^{~N}))S(\Lambda_{L}^{~P})(X_{(b)Q}\otimes 
X_{(a)N})
\end{eqnarray*}
where $a,b = \pm$ or $\mp$. From the later equation we deduce
\begin{eqnarray*}
\sigma (X_{(a)L} \otimes X_{(b)K})=F_{(b)K}^{~~N}
(S(\Lambda_{L}^{~M}))(X_{(b)N}\otimes X_{(a)M}).
\end{eqnarray*}
For the spinors, it suffices to replace $\Lambda_{L}^{~K}$ by
$M_{\alpha}^{~\beta}$ or  $M_{\dot{\alpha}}^{~\dot{\beta}}$ and $X_{(a)I}$
by $\theta_{(a)\alpha}$ or $\theta_{(a)\dot{\alpha}}$.
The symmetrization of the product is defined as
\begin{eqnarray}
X_{(a)L}X_{(b)K} &=& F_{(b)K}^{~~N}(S(\Lambda_{L}^{~M}))(X_{(b)N}X_{(a)M}),
\nonumber\\
\theta_{(a)\alpha}X_{(b) K}=F_{(b) K}^{~~N}(S(M_{\alpha}^{~\beta}))(X_{(b)N}
\theta_{(a)\beta})~&or&~
\theta_{(a)\dot{\alpha}}X_{(b) K}=F_{(b) K}^{~~N}
(S(M_{\dot{\alpha}}^{~\dot{\beta}}))(X_{(b)N}\theta_{(a)\dot{\beta}}).
\end{eqnarray}
From this, we state
\vskip0.25truecm
{\bf Theorem (4,1):}\\
$The~norm~ of~ X_{\pm},~G_{\pm}^{~IJ}X_{\pm I}X_{\pm J}~ is~ invariant~ and~
central.$\\
\vskip0.25truecm
$Proof.$: By construction, $X_{\pm I}$ is right invariant. As a consequence
of the orthogonality condition of the generators of $\cal L$, we may easily
see from the transformations of
$X_{\pm I}$ (34) that $G_{\pm}^{~IJ}X_{\pm I}X_{\pm J}$ is left invariant.
From (47) and (48), we get, for any $a \in {\cal A}$,
\begin{eqnarray}G_{\pm}^{~IJ}X_{\pm I}X_{\pm J}a =
(a \star G_{\pm}^{~IJ}(F_{\pm J}^{~~L} \star F_{\pm I}^{~~K}))
X_{\pm K}X_{\pm L}=\nonumber\\
(a \star \varepsilon)G_{\pm}^{~KL}X_{\pm K}X_{\pm L} = aG_{\pm}^{~KL}X_{\pm K}
X_{\pm L}
\end{eqnarray}
which shows that the norm commutes with any $a \in \cal A$.\\
From (48) and (49), we obain
\begin{eqnarray}
X_{(a) P}G_{\pm}^{~IJ}X_{\pm I}X_{\pm J} = G_{\pm}^{~IJ}F_{\pm I}^{~~N}
(S(\Lambda_{P}^{~M})) F_{\pm J}^{~~K}(S(\Lambda_{M}^{~L}))X_{\pm N}X_{\pm K}
X_{(a) L}=\nonumber\\
G_{\pm}^{~IJ}(F_{\pm J}^{~~K} \star F_{\pm I}^{~~N})(S(\Lambda_{P}^{~L}))
X_{\pm N}X_{\pm K}X_{(a) L}= G_{\pm}^{~NK}X_{\pm N}X_{\pm K}X_{(a)P}.
\end{eqnarray}
The same results may be obtained by replacing $X_{(a) P}$ by $\theta_{(a)\alpha}$ 
or $\theta_{(a)\dot{\alpha}}$ which show that the norm also commutes with the
quantum coordinates and the spinors.  Q.E.D.\\
\vskip0.25truecm
{$Remark (4,1):$}\\
- $G_{\pm}^{~LK}X_{\pm L}X_{\pm K}$ is biinvariant and real. Since it 
commutes with everything, it is of the form 
$\lambda I_{\cal L}$ with $\lambda$ is a real number.\\
- The quantities $G_{\pm}^{~IJ}X_{\mp I}X_{\mp J}$ and 
$G_{\pm}^{~IJ}X_{\mp I}X_{\pm J}$ are biinvariant but not central. 
\vskip0.25truecm
{\bf Acknowledgments.} I am grateful to M. Dubois-Violette for his kind
interest and hepful suggestion. I also thank N. Touhami for valuable
discussions and S. Balaska for his help which allowed me to have results
from the computer MAPLE program.\\
\vskip0.25truecm
{\bf References:}\\
1)A. Connes, Publ, IHES 62(1986)257, Noncommutative Geometry, Academic Press,
inc. (1994).\\
2)R. Coquereaux, G. Esposito-Farese and G. Vaillant, Nucl. Phys. B353(1991)
689.\\
3)M. Dubois-Violette, R. Kerner and J. Madore, Class. Quant. Grav.
6(1989)1709, J. Math. Phys. 31(1990)323.\\
4)S. L. Woronowicz, Commun. Math. Phys. 122(1989)125.\\
5)P. Podles and S.L. Woronowicz, Commun. Math. Phys. 130(1990)381.\\
6)U. Carow-Watamura, M. Schlieker, M. Scholl and S. Watamura, Z. Phys. C-
$partcles~and~fields$ 48(1990)159, Int. J. Mod. Phys. A6(1991)3081.\\
7)O. Ogiesvestsky, W. B. Schmidke, J. Wess and B. Zumino, Commun. Math.
Phys. 150(1992)495.\\
8)P. Podles and S. L. Woronowicz, Commun. Math. Phys. 178(1996)61,
Commun. Math. Phys. 185(1997)325.\\
9)M. Dubois-Violette and G. Launer, Phys. Lett. B245(1990)175.\\
10)M. Lagraa, Int. J. Mod. Phys. A11(1996)699.\\
11)A. B. Hammou and M. Lagraa, J. Math. Phys. 38(1997)4462.
\end{document}